\documentclass[10pt,a4j]{article}

\usepackage{lamatome}

\begin{document}
\title{Construction of schemes
over $\FF_{1}$, and over idempotent semirings:
towards tropical geometry}
\author{Satoshi Takagi}
\maketitle

\begin{abstract}
In this paper,
we give some categorical description
of the general spectrum functor,
defining it as an adjoint
of a global section functor.
The general spectrum functor
includes that of $\FF_{1}$ and
of semirings.
\end{abstract}

\tableofcontents

\setcounter{section}{-1}
\section{Introduction}
The goal of this paper
is to provide a explicit
description of the
mechanism of the spectrum functor;
how we construct a topological
space and a structure sheaf
from an algebraic object,
and why it behaves so nicely.

At first, we were only planning
to give some definitions of
schemes constructed from 
idempotent semirings 
motivated from the tropical geometry
(cf. \cite{IME}),
since still there is no concrete
definition of higher dimensional
tropical varieties:
in \cite{Mik}, Mikhalkin
is giving a definition,
but is very complicated
and not used in full strength yet,
so we cannot say that the definition
is appropriate yet.
Also, the algebraic theories of
semirings are still inadequate,
compared to those of normal algebras,
though there are already many demands
from various workfields. (cf. \cite{Lit}).
Even the tensor products has been
constructed only recently (cf. \cite{LMS}).
So, we had to start from the
beginning: giving a foundation
of idempotent algebras.

Later, we noticed that there should be a comprehensive
theory including  the theory of schemes
over $\FF_{1}$, of which field is now
in fashion (cf. \cite{TV}, ~\cite{CC}, ~\cite{PL}, etc.).

We chose the definition of 
commutative idempotent monoid to be infinitely additive.
There are two reasons.
First, if we restrict
our attention to finite additivity,
then the existence of all
adjoint functors are already assured,
by the general theory of algebraic systems.
Secondly,
we can handle topological spaces
in this framework:
a sober space can be regarded
as a semiring with idempotent multiplication.
Also, we can make it clear how we obtain
a infinite operation (i.e.
the intersection of closed subsets)
from algebras with only finitary operators.

This fact shows that,
there are adjoints between
the three categories:
that of algebraic complete idealic semirings,
that of algebraic idealic schemes, 
and that of algebraic sober spaces:
\[
\cat{alg.IRng$^{\dagger}$} \leftrightarrows 
\cat{alg.ISch}^{\op} \leftrightarrows 
\cat{alg.Sob}^{\op}.
\]
This indicates that,
the genuine theme of
algebraic geometry and
arithmetic geometry 
(including schemes over $\FF_{1}$) is
no longer the investigation
of the link between the category
of algebras and the category of
topological spaces,
but that of the link between
the category of algebras
and the category of idempotent semirings.
Roughly speaking,
the whole workfield may be thought as
contained in the tropical geometry,
if you take the tropical geometry
as a study of semirings.

After a while, we realized that
there are strong resemblance
between the theory of lattices and
the argument in the first half of this article.
We must admit that there are no new ideas
contained in this article;
most of the techniques and
theorems have been already established
decades ago. 
Therefore, readers who attempt
to study in this workfield is strongly
recommended to overlook the lattice theory
(for example, ~\cite{MMT}, \cite{JH}).
However,
we decided to contain definitions
and propositions as self-contained as possible,
since there seems to be some
discrepancy of languages between
each workfields,
and one of the purpose
of this paper is to introduce connections
between the lattice theory,
tropical geometry,
and schemes over $\FF_{1}$.

We were irresolute for a 
while, of which definition of schemes
would be the most general and natural one.
For this,  we decided to
give the weakest (in other words,
the most general) definition in which
the global section
functor $\Gamma$ becomes an adjoint.
This clarifies the condition
of which algebra could give a
topological space endowed with
a structure sheaf.

Also, we know that we can further
extend the arguments and constructions
in this article by replacing commutative monoids
with symmetric \linebreak
monoidal categories,
but we dared not,
for it would be too abstract for
the present demands.

\S 1 is a summary of algebras
and complete algebras.
The definition and arguments are
all standard -- in other words, preliminary -- in the lattice theory.
In \S 2, we deal with $R$-modules,
where $R$ is a complete semiring.
The results cannot be deduced from
\S 1, so we had to argue seperately.
In \S 3, we focus on congruence relations,
which is important when considering
localizations.
In \S 4 we deal with topological spaces,
defining the spectrum functor by an adjoint.
This is also standard in the lattice theory.
\S 5 may be the most valuable part,
in which we define $\scr{A}$-schemes,
and define the spectrum functor
(not the same one of \S 4, but
also endows the structure sheaf to
the underlying space) as an adjoint.
All the adjoints we have constructed are illustrated
in the end of 5.3.

\textit{Acknowledgements.}
We would like to thank
Professor Moriwaki for giving me chances
to work on this subject.
Also, the author is
grateful to Doctor Ikoma,
for giving useful advises and
encouragements.

\subsection{Notation and conventions}

Throughout this paper,
all monoids are assumed associative,
unital, commutative.

We fix a universe, and
all sets and algebraic objects are elements
of this universe.
For any set $X$, we denote the
power set of $X$ by $\scr{P}(X)$.

We frequently make use of the
notation $\sum^{<\infty}$:
this means that given a infinite sum,
there exists a finite subindex 
which preserve the equality and inequality.
For example,
\[
x \leq  \sum_{\lambda} a_{\lambda}
\Rightarrow x \leq \sum^{<\infty}_{\lambda} a_{\lambda}
\]
means that, if $x<\sum_{\lambda} a_{\lambda}$,
then there are finitely many $\lambda_{1},\cdots, \lambda_{n}$
such that $x <\sum_{i=1}^{n} \lambda_{i}$.

We also make use of the notation $\cup^{<\infty}$.

\subsection{Preliminaries}
\begin{enumerate}
\item
We frequently make use of the categorical languages in ~\cite{CWM},
especially that of adjoint functors:
Let $U: \scr{A} \to \scr{C}$,
$F: \scr{C} \to \scr{A}$ be two functors
between two categories $\scr{A}$ and $\scr{C}$.
We say $F$ and $U$ are \textit{adjoint},
or $F$ is the \textit{left adjoint of $U$},
if there is a natural isomorphism of functors
\[
\Hom_{\scr{C}}(c,Ua) \simeq \Hom_{\scr{A}}(Fc,a):
\scr{C}^{\op} \times \scr{A} \to \cat{Set}.
\]
This is equivalent to saying that
there are natural transformations
$\epsilon:\Id_{\scr{C}} \Rightarrow  UF$ (the \textit{unit}),
and $\eta:FU \Rightarrow \Id_{\scr{A}}$ (the \textit{counit})
which satisfy
$\eta F \circ F\epsilon =\Id_{F}$ and
$U \eta \circ \epsilon U=\Id_{U}$:
\[
\xymatrix{
F \ar[d]_{F\epsilon} \ar@{=}[dr] & \\
FUF \ar[r]_{\eta F} & F
}
\xymatrix{
U \ar[r]^{\epsilon U} \ar@{=}[rd] & UFU \ar[d]^{U\eta} \\
& U
}
\]
\item
We will briefly overview
the general theory of algebraic varieties
(in the lattice theoritic sense).
Here we also use the language of \cite{CWM}.
An \textit{algebraic type} $\sigma$
is a pair $\langle \Omega, E \rangle$,
where $\Omega$ is a set of finitary operators
and $E$ is a set of identities.
Let $D$ be the set of \textit{derived operators},
i.e. the minimal set of operators including $\Omega$,
and closed under compositions and substitutions.
A $\sigma$-algebra is a set
$A$, with an action of $\Omega$ on $A$,
satisfying the identities.
A \textit{homomorphism of $\sigma$-algebras}
is a map $f:A \to B$ between two
$\sigma$-algebras, preserving the actions
of $\Omega$.
We denote by $\cat{$\sigma $-alg}$
the category of small $\sigma$-algebras.
Here are some facts about algebraic types:
\begin{enumerate}
\item The category $\cat{$\sigma $-alg}$
is small complete and small co-complete.
\item Let $f:\sigma \to \tau$
be a \textit{homomorphism of algebraic types},
i.e. $f$ is a map $\Omega_{\sigma} \to D_{\tau}$
of operators which satisfies $f(E_{\sigma}) \subset E_{\tau}$.
We say that $\tau$ is stronger than $\sigma$
in this case, and denote by $\sigma \leq \tau$ for brevity.
Then, the underlying functor
$U:\cat{$\tau $-alg} \to \cat{$\sigma $-alg}$
has a left adjoint.
In particular, there is a functor
$F:\cat{Set} \to \cat{$\sigma $-alg}$
for any algebraic type $\sigma$,
and the $\sigma$-algebra $F(S)$
is called the \textit{free $\sigma$-algebra}
generated by $S$.
\end{enumerate}

\end{enumerate}

\section{Complete algebras}

Almost all the results
are already established in the lattice
theory. Here, we are just translating
it into a somewhat algebro-geometric
language. The reader
who is well familiar with the lattice theory
may skip.
\subsection{Complete sets}
\begin{Def}
\begin{enumerate}
\item A \textit{multi-subset}
of a set $X$ is simply, a map $f: \Lambda \to X$
from any set $\Lambda$.
We denote by $\bsP(X)$
 the set of multi-subset of $X$
modulo isomorphism.
Often, a multi-subset is identified
with its image.
\item A map $\sup: \bsP(X) \to X$
is a \textit{supremum map} if it satisfies:
\begin{enumerate}
\item The map factors through the power set of $X$,
namely, there is a map $\scr{P}(X) \to X$
making the following diagram commutative:
\[
\xymatrix{
\bsP(X) \ar[r]^{\sup} \ar[d]_{\Imag} & X \\
\scr{P}(X) \ar[ur]
}
\]
where $\Imag$ sends a map $\Lambda \to X$
to its image.
\item It is (infinitely) idempotent:
It sends a constant map (from a non-empty set)
to the unique element in its image.
\item It is (infinitely) associative:
If $f:\Lambda \to X$ and 
$f_{i}:S_{i}\to X$ are multi-subsets of $X$
satisfying $\Imag f=\cup_{i} \Imag f_{i}$,
then $\sup f=\sup \{\sup f_{i}\}_{i}$.
\end{enumerate}
We often write $x \oplus y$ instead of $\sup\{x,y\}$.
Note that when given a (not neccesarily infinite)
associative commutative idempotent operator
$\oplus$ on $X$,
we can define a preorder on $X$,
by defining
\[
a \leq b \Leftrightarrow a \oplus b=b.
\]
Conversely, the supremum map
gives the supremum with respect to this preorder.
The maximal element $1$ of $X$
is the \textit{absorbing element} of $X$,
i.e. $\sup (S \cup \{1\})=1$ for any subset $S$ of $X$.
Also, there is an infimum
of any subset $S$ of $X$:
\[
\inf S=\sup \{ x \in X \mid \text{$x \leq s$ for any $s \in S$}\},
\]
and the minimal element $0$ of $X$
is the \textit{unit}, i.e. $\sup S \cup \{0\}=\sup S$
for any subset $S$ of $X$.
\item
A \textit{complete idempotent monoid}
is a set endowed with a supremum map.
\item Let $A,B$  be two complete idempotent
monoids. A map $f:A \to B$ is
a \textit{homomorphism of complete idempotent
monoids} if it commutes with the supremum map:
\[
\xymatrix{
\bsP(A) \ar[r]^{f_{*}} \ar[d]_{\sup} & \bsP(B) \ar[d]^{\sup} \\
A \ar[r]^{f} & B
}
\]
and sending $1$ to $1$.
Note that when this holds,
$f$ sends $0$ to $0$,
since $0=\sup \emptyset$.
We denote by $\cat{CIM}$
the category of complete idempotent monoids.
\end{enumerate}
\end{Def}
The notion "complete idempotent
monoid" is used in ~\cite{LMS},
but in the lattice theory, this is
refered to as a \textit{complete lattice}.
We will use the former notation.

\begin{Def}
The algebraic type $\sigma$
with a binomial operator $+$ is
\textit{pre-complete} if:
\begin{enumerate}
\item 
$+$ is an associative, commutative idempotent
operator with a unit $0$ and absorbing element $1$,
\item
There are infimums for any two elements $a,b$:
$\inf(a,b) \leq a,b$, and if $x \leq a,b$
then $x \leq \inf(a,b)$.
\item If $\phi$ is another $n$-ary
operator, then it is $n$-linear with respect to $+$:
\begin{multline*}
\phi(x_{1},\cdots, x_{i}+x_{i}^{\prime},\cdots, x_{n}) \\
=\phi(x_{1},\cdots, x_{i},\cdots, x_{n})
+\phi(x_{1},\cdots, x_{i}^{\prime},\cdots, x_{n}).
\end{multline*}
\end{enumerate}
\end{Def}
The condition (1) and (2)
are equivalent to saying that a
$\sigma$-algebra is a \textit{distributive lattice}
with respect to $+$ and $\inf$.

Note that any complete idempotent
monoid is already pre-complete,
if we restrict the supremum map
to finite subsets.
Also, see that $\inf$ is defined algebraically,
i.e.
\begin{enumerate}[(a)]
\item $a+\inf(a,b)=a$.
\item $\inf(x+a,x+b)+x=\inf(x+a,x+b)$.
\end{enumerate}

\begin{Def}
Let $\sigma$ be a pre-complete
algebraic type with respect to $+$.
\begin{enumerate}
\item We denote by $\cat{$\sigma $-alg}$
the category of $\sigma$-algebras.
\item A \textit{complete $\sigma$-algebra} is
a $\sigma$-algebra $A$ satisfying:
\begin{enumerate}
\item $A$ is a complete idempotent monoid,
and the supremum map coincides with $+$,
when restricted to any finite subset of $A$.
\item Any $n$-ary operator $\phi:A^{n} \to A$
is $n$-linear with respect to $+$:
\[
\phi(x_{1},\cdots, \Sum_{j}x_{ij},\cdots, x_{n})
=\Sum_{j}\phi(x_{1},\cdots, x_{ij},\cdots, x_{n}).
\]
\end{enumerate}
\item A \textit{homomorphism $f:A \to B$
of complete $\sigma$-algebras}
is a homomorphism of $\sigma$-algebras
and homomorphism of complete idempotent monoids.
\item We denote by $\cat{$\sigma^{\dagger} $-alg}$
the category of complete $\sigma$-algebras.
We refer to $\sigma^{\dagger}$ as a
\textit{complete algebraic type}.
\end{enumerate}
\end{Def}

We will see that
the underlying functor
$U: \cat{$\sigma^{\dagger} $-alg} \to \cat{$\sigma $-alg}$
has a left adjoint,
but we will analyse this left adjoint more
precisely, for future references.

\begin{Def}
Let $A$ be a complete $\sigma$-algebra.
\begin{enumerate}
\item Let $a$ be an element of $A$.
A subset $S$ of $A$ is a \textit{covering of $a$},
if $a \leq \sup S$.
\item An element $a$ of $A$
is \textit{compact}, if any covering
of $a$ has a finite subcovering of $a$.
\item We say that $A$ is \textit{algebraic},
if the following holds:
\begin{enumerate}
\item For any element $a$ of $A$ is \textit{algebraic},
i.e.
$a$ has a covering $S$ which consists
of compact elements.
\item Any operator (including the infimum of finite elements)
$\phi:A^{n} \to A$
preserves compactness, i.e. 
$\phi(x_{1},\cdots,x_{n})$ is compact
if $x_{i}$'s are compact elements.
\end{enumerate}
\item A \textit{homomorphism
$f:A \to B$ of algebraic complete $\sigma$-algebras}
is a homomorphism of complete $\sigma$-algebras,
sending any compact element to a compact element.
\item We denote by $\cat{alg.$\sigma^{\dagger}$-alg}$
the category of algebraic complete $\sigma$-algebras.
\end{enumerate}
\end{Def}
These notation come
from the lattice theory.

Note that any constant (regarded as a $0$-ary operator)
in an algebraic complete $\sigma$-algebra is compact,
by definition. In particular,
the absorbing element $1$ is compact.

\begin{Prop}
\label{prop:adj:completion}
Let $\sigma$ be a pre-complete algebraic type.
Then, the underlying functor
$U: \cat{$\sigma^{\dagger} $-alg} \to \cat{$\sigma $-alg}$
has a left adjoint $\comp$.
Further, $\comp$ factors through
$\cat{alg.$\sigma^{\dagger}$-alg}$:
\[
\xymatrix{
\cat{$\sigma$-alg} \ar[r]^{\comp} \ar[d]_{\comp^{\prime}} 
	& \cat{$\sigma^{\dagger}$-alg} \\
\cat{alg.$\sigma$-alg} \ar[ru]_{U^{\prime}} 
}
\]
\end{Prop}

\begin{proof}
Let $A$ be a $\sigma$-algebra.
A \textit{filter} $F$ on $A$ is
a non-empty subset of $A$ satisfying:
\[
x,y \in F \Leftrightarrow x+y \in F.
\]
We denote by $\langle x \rangle$
the filter generated by an element $x \in A$.
Let $A^{\dagger}$ be the set of all filters on $A$.
Given a family of filters $\scr{F}=\{F_{\lambda}\}_{\lambda}$,
the supremum of $\scr{F}$ is the
filter generated by $F_{\lambda}$'s,
i.e. $\xi \in \sup\scr{F}$ if and only if
there are \textit{finite} number of $x_{\lambda} \in F_{\lambda}$'s
such that $\sum_{\lambda} x_{\lambda} \geq \xi$.
Then, the unit is $\{0\}$ and the absorbing
element is $A$.
We can easily see that that $1$
is compact, and any element is algebraic:
a compact filter is precisely, a finitely generated filter.
Let $\phi:A^{n} \to A$ be another operator of $A$.
The operator $\phi^{\dagger}:(A^{\dagger})^{n} \to A^{\dagger}$
associated to $\phi$ is defined by
\[
\phi^{\dagger}(F_{1},\cdots,F_{n})
=\sum_{x_{i} \in F_{i}}\langle \phi(x_{1},\cdots, x_{n}) \rangle.
\]
This sends a $n$-uple of compact filters
to a compact filter.
This gives a structure of an algebraic complete $\sigma$-algebra on $A^{\dagger}$.
Given a homomorphism $f:A \to B$
of $\sigma$-algebras,
a homomorphism $f^{\dagger}:A^{\dagger} \to B^{\dagger}$
of complete $\sigma$-algebras is given by
\[
f^{\dagger}(F)=\sum_{a \in F} \langle f(a) \rangle.
\]
It is easy to see that $f^{\dagger}$ sends
a compact filter to a compact filter.
Hence, we have a functor
$\comp^{\prime}:\cat{$\sigma $-alg} \to 
\cat{alg.$\sigma^{\dagger}$-alg}$.
Set $\comp=U^{\prime} \circ \comp^{\prime}$,
where $U^{\prime}: \cat{$\sigma^{\dagger}$-alg}
\to \cat{alg.$\sigma^{\dagger}$-alg}$
is the underlying functor.

Finally, we will show that $\comp$ is the left adjoint
of $U$:
the unit $\epsilon:\Id_{\cat{$\sigma $-alg}} 
\Rightarrow U\circ \comp$
is given by $A \ni a \mapsto \langle a \rangle \in A^{\dagger}$.
$\eta:\comp \circ U \Rightarrow
 \Id_{\cat{alg.$\sigma^{\dagger} $-alg}}$
is given by $B^{\dagger} \ni F \mapsto \sup F \in B$.
\end{proof}
\begin{Rmk}
The functor $\comp^{\prime}$ constructed above is
\textit{not} the left adjoint of the underlying functor
$\cat{alg.$\sigma^{\dagger}$-alg} \to \cat{$\sigma $-alg}$:
the counit $\eta$ is not algebraic, i.e.
it does not necessarily preserve compactness.
However,
it is an equivalence of categories: see below.
\end{Rmk}

\begin{Prop}
\label{prop:comp:isom:cpt}
Let $\sigma$ be a pre-complete algebraic type.
Then the above functor
$\comp^{\prime}$ gives an equivalence between
the category
of $\sigma$-algebras and
the category of algebraic complete $\sigma$-algebras.
\end{Prop}
\begin{proof}
We will construct a functor
$U_{\cpt}:\cat{alg.$\sigma^{\dagger}$-alg} \to \cat{$\sigma $-alg}$
as follows:
for an algebraic complete $\sigma$-algebra
$R$, let $R_{\cpt}$ be the set
of compact elements of $R$.
This set has the natural induced structure
of a $\sigma$-algebra, since all the operators
are algebraic.
Also, given a homomorphism $f:A \to B$
of algebraic complete $\sigma$-algebras,
we obtain a homomorphism
$f_{\cpt}:A_{\cpt} \to B_{\cpt}$
of $\sigma$-algebras.
Hence, sending $R$ to $R_{\cpt}$ gives a functor 
$U_{\cpt}:\cat{alg.$\sigma $-alg} \to \cat{$\sigma $-alg}$.

We will see that $U_{\cpt}$ is the
inverse of $\comp^{\prime}$.
The unit $\epsilon: \Id_{\cat{$\sigma $-alg}} \Rightarrow
 U_{\cpt} \circ \comp^{\prime}$
is given by
\[
A \ni a \mapsto \langle a \rangle \in (A^{\dagger})_{\cpt}.
\]
This is an isomorphism:
the inverse is given by $F \mapsto \sup F$.
Note that this is well defined,
since a compact filter is finitely generated,
hence only finitely many elements is
involved when taking its supremum.

The counit
$\eta: \comp^{\prime} \circ U_{\cpt} 
\Rightarrow \Id_{\cat{alg.$\sigma^{\dagger}$-alg}}$
is given by
\[
(B_{\cpt})^{\dagger} \ni F \mapsto \sup F \in B.
\]
This is well defined, since all the filters
are generated by compact elements of $B$,
hence $\eta$ is algebraic.
The inverse of $\eta$ is given by
\[
B \ni b \mapsto \sum_{b^{\prime}\leq b} \langle b^{\prime} \rangle
\in (B_{\cpt})^{\dagger},
\]
where $b^{\prime}$ runs through
all the compact
elements smaller than $b$.
It is clear that $\eta^{-1}$ preserves compactness.
Also, $\eta^{-1}$ preserves any operator $\phi$ of $\sigma$,
since $\phi$ preserves compactness.
Hence, $\eta^{-1}$ is well defined as
a homomorphism of algebraic complete $\sigma$-algebras.
\end{proof}
\begin{Cor}
\label{cor:alg:sigma:alg:comp}
Let $\sigma$ be a pre-complete algebraic type.
\begin{enumerate}
\item 
The underlying functor
$U^{\prime}: \cat{alg.$\sigma^{\dagger}$-alg} \to
\cat{$\sigma^{\dagger}$-alg}$ has a right adjoint:
it is $\alg=\comp^{\prime} \circ U$,
where $U: \cat{$\sigma^{\dagger}$-alg} \to \cat{$\sigma $-alg}$
is the underlying functor. 
\item The category $\cat{alg. $\sigma^{\dagger}$-alg}$
is small complete and small co-complete.
\item Let $\tau$ be another pre-complete algebraic type
stronger than $\sigma$.
Then we have a following natural commutative
diagram of functors:
\[
\xymatrix{
\cat{$\tau $-alg} \ar[r]^>(0.7){\simeq} \ar[d]
& \cat{alg.$\tau^{\dagger}$-alg} \ar[d] \\
\cat{$\sigma $-alg} \ar[r]_>(0.7){\simeq} & \cat{alg.$\sigma^{\dagger}$-alg} \\
}
\]
where the downward arrows are the underlying functors.
Hence, the underlying functor
$U:\cat{alg.$\tau^{\dagger}$-alg} \to
\cat{alg.$\sigma^{\dagger}$-alg}$
has a left adjoint.
\end{enumerate}
\end{Cor}

\begin{Prop}
\label{prop:comp:sigma:complete}
Let $\sigma$ be a pre-complete algebraic type.
\begin{enumerate}
\item
The category $\cat{$\sigma^{\dagger} $-alg}$
is small complete.
\item
The category $\cat{$\sigma^{\dagger} $-alg}$
is small co-complete.
\end{enumerate}
\end{Prop}
\begin{proof}
\begin{enumerate}
\item
We only need to verify
that small products and equalizers exist,
but this is just the analogue of the case
of usual algebras.
\item 
We need to verify
that small co-products and  co-equalizers exist.
First, we will verify the existence
of co-equalizers.
Let $A \stackrel{f,g}{\rightrightarrows} B$
be two homomorphisms between two
$\sigma$-algebras $A$ and $B$.
Let $\mathfrak{a}$ be the \textit{congruence relation}
generated by $f$ and $g$, i.e.
the minimal equivalence relation
on $B$ satisfying:
\begin{enumerate}
\item $(f(a),g(a)) \in \mathfrak{a}$ for any $a \in A$.
\item If $(a_{i},b_{i}) \in \mathfrak{a}$ for all $i$,
then $(\phi(a_{1},\cdots,a_{n}), \phi(b_{1},\cdots,b_{n}))
\in \mathfrak{a}$ for any $n$-ary operator $\phi$.
\item If $(a_{\lambda},b_{\lambda}) \in \mathfrak{a}$
for all $\lambda$, then
$(\sum a_{\lambda},\sum b_{\lambda}) \in \mathfrak{a}$.
\end{enumerate}
Then, it is clear that $B/\mathfrak{a}$
is the co-equalizer of $f$ and $g$.

Next, we will show that small co-products exist.
Let $\scr{A}=\{A_{\lambda}\}_{\lambda}$
be a small family of complete $\sigma$-algebras.
Let $S$ be the epic undercategory of $\scr{A}$,
i.e. its objects are pairs
$\langle B, \{f_{\lambda}:A_{\lambda} \to B\}\rangle$
satisfying:
\begin{enumerate}
\item $B$ is a complete $\sigma$-algebra.
\item $f_{\lambda}$ is a homomorphism of complete
$\sigma$-algebras.
\item The image of $f_{\lambda}$'s
generate $B$ as a complete $\sigma$-algebra.
\end{enumerate}
Then, we see that the set $S$ is small,
and complete when regarded as a category,
since $\cat{$\sigma^{\dagger}$-alg}$
is complete. The coproduct $\amalg A_{\lambda}$
is then defined by the limit of $S$.
\end{enumerate}
\end{proof}

\begin{Prop}
\label{prop:adj:comp:strong}
Let $\sigma$, $\tau$
be two pre-complete algebraic types,
and assume $\tau$ is stronger than $\sigma$.
Then, the underlying functor
$U_{v}: \cat{$\tau^{\dagger}$-alg} \to \cat{$\sigma^{\dagger}$-alg}$
has a left adjoint.
\end{Prop}
\begin{proof}
Let
$F_{\sigma}:\cat{$\sigma $-alg} \leftrightarrows  
\cat{$\sigma^{\dagger}$-alg}:U_{\sigma}$,
$F_{\tau}:\cat{$\tau $-alg} \leftrightarrows  
\cat{$\tau^{\dagger}$-alg}:U_{\tau}$,
and 
$F_{\sigma}:\cat{$\sigma $-alg} \leftrightarrows  
\cat{$\tau $-alg}:U_{\sigma}$
be adjoints, respectively:
\[
\xymatrix{
\cat{$\sigma $-alg} \ar@{<->}[r]^{u} \ar@{<->}[d]^{\sigma} & 
\cat{$\tau $-alg} \ar@{<->}[d]^{\tau} \\
\cat{$\sigma^{\dagger}$-alg} \ar@{<.>}[r]^{v} &
\cat{$\tau^{\dagger}$-alg}
}
\]
Given a complete $\sigma$-algebra $A$,
set $A^{\prime}=F_{u}U_{\sigma}A \in \cat{$\tau $-alg}$.
A filter $F$ on $A^{\prime}$ is \textit{infinite}
with respect to $A$, if
$x_{\lambda} \in F \cap A$ implies $\sum x_{\lambda} \in F$.
Let $F_{v}(A)$ be the set of infinite
filters on $A^{\prime}$.
When given a family $\{F_{\lambda}\}$
of infinite filters, the supremum filter
is the infinite filter generated by $F_{\lambda}$'s.
The rest of the construction of
the complete $\sigma$-algebra structure on $F_{v}(A)$
is analogous to Proposition ~\ref{prop:adj:completion}.
Hence, we have a functor
$F_{v}:\cat{$\sigma^{\dagger}$-alg} 
\to \cat{$\tau^{\dagger}$-alg}$.
We will show that this is the left adjoint of $U_{v}$.
The unit $\epsilon:\Id_{\cat{$\sigma^{\dagger}$-alg}} 
\Rightarrow U_{v}F_{v}$
is given by $A \ni a \mapsto \langle a \rangle \in F_{v}(A)$.
Note that this preserves the supremum map,
from the advantage of using infinite filters.
The counit $\eta: F_{v}U_{v} \Rightarrow \Id_{\cat{$\tau^{\dagger}$-alg}}$
is given by $F \mapsto \sup F$.
\end{proof}

\subsection{Semirings}
\begin{Def}
\begin{enumerate}
\item An algebraic system $(R,+,\times)$
is a \textit{semiring} if:
\begin{enumerate}
\item $R$ is a monoid with respect to $+$ and $\times$.
\item The \textit{distribution law} holds:
\[
(a+b)\cdot c=(a\cdot c)+(b \cdot c).
\]
\end{enumerate}
\item A semiring $R$ is \textit{pre-complete},
if it is pre-complete with respect to $+$,
and if the multiplicative unit is the absorbing element
with respect to $+$.
\item We denote by $\cat{SRng}$ (resp. $\cat{PSRng}$),
the category of semirings (resp. pre-complete semirings).
\end{enumerate}
\end{Def}

\begin{Def}
\begin{enumerate}
\item The initial object of $\cat{PSRng}$
consists of two elements: $1$ and $0$.
We refer to this semiring as $\FF_{1}$,
the \textit{field with one element}.
Here in the category of semirings,
$\FF_{1}$ is no longer a makeshift concept.

\item The terminal object of $\cat{SRng}$
exists, and consists of a unique element $1=0$.
We refer to this semiring as the \textit{zero semiring}.
\item Let $R$ be a complete semiring.
Let
$U_{m}:\cat{$R$/SRng$^{\dagger}$} \to \cat{Mnd}$
be the underlying functor, sending
a complete $R$-algebra to its multiplicative monoid.
Then, $U_{m}$ has a left adjoint $F_{m}$.
For a monoid $M$, $R[M]:=F_{m}(M)$
is the \textit{complete monoid semiring}
with coefficient $R$.
In particular,
when $M$ is the free monoid generated by a set $S$,
$R[x_{s}]_{s\in S}=F(M)$
is the \textit{complete polynomial semiring}
with coefficient $R$. 
\item Let $A$, $B$, $R$ be complete semirings,
with homomorphisms $R \to A$ and $R \to B$.
Then we can define the \textit{tensor product
of $A$  and $B$ over $R$} which we denote
by $A \otimes_{R} B$.
Of course, we have infinite tensor products too.
\end{enumerate}
\end{Def}

\begin{Rmk}
The constructions of various
algebras are also valid in the algebraic complete
case. The arguments and notations will be
just the repetition, so we will skip it.
\end{Rmk}

\section{Modules over a semiring}
We will construct
a general fundamental theorems
on $R$-modules, where $R$
is a complete semiring.
Note that $\cat{$R$-mod}$ is a little
beyond from what we have been calling
complete algebras:
the definition of the scalar action is not
algebraic, since we must consider infinite
sums of elements of the coefficient ring $R$.
However, the arguments are analogous.
We just modify the definition of the filters when necessary.

\subsection{$R$-modules}
\begin{Def}
Let $R$ be a complete semiring.
\begin{enumerate}
\item An \textit{$R$-module} $M$ is
a complete idempotent monoid endowed
with a scalar operation $R \times M \to M$
which satisfies the following:
\begin{enumerate}
\item $1\cdot x=x$, $0 \cdot x=0$,
$a \cdot 0=0$ for any $x \in M$ and $a \in R$.
\item The \textit{distribution law} holds:

$a(\sum_{\lambda}x_{\lambda})=\sum_{\lambda}(ax_{\lambda})$
for any $a \in R$, $x_{\lambda} \in M$.

$(\sum_{\lambda}a_{\lambda})x=\sum_{\lambda}(a_{\lambda}x)$
for any $a_{\lambda} \in R$, $x \in M$.
\item $(ab)x=a(bx)$ for any $a,b \in R$ and $x \in M$.
\end{enumerate}
\item A map $f:M \to N$ between two $R$-modules
is a \textit{homomorphism of $R$-modules}
if it is a homomorphism of complete idempotent
monoids, and preserves the scalar operation.
\item We denote by $\cat{$R$-mod}$
the category of $R$-modules.
\end{enumerate}
\end{Def}

\begin{Prop}
\begin{enumerate}
\item The category of $\FF_{1}$-modules coincides
with the category of complete idempotent monoids.
\item Let $R$ be a complete semiring.
For any two $R$-modules $M$ and $N$,
$\Hom_{\cat{$R$-mod}}(M,N)$ has the
natural structure of a $R$-module.
\end{enumerate}
\end{Prop}
The proofs are obvious.

\begin{Prop}
Let $R$ be a complete semiring.
\begin{enumerate}
\item The category $\cat{$R$-mod}$ is complete.
\item The category $\cat{$R$-mod}$ is co-complete.
\end{enumerate}
\end{Prop}
\begin{proof}
\begin{enumerate}
\item It is just the same argument
of classical modules.
\item It suffices to show
that there is small co-products and coequalizers.
The construction of a coequalizer is easy.
We will construct a co-product
for a small family $\{M_{\lambda}\}$ of $R$-modules.
Set $M_{\infty}=\prod M_{\lambda}$:
the set theoritic product of $M_{\lambda}$'s.
We will show that this is also the co-product.
The inclusion $M_{\lambda} \to M_{\infty}$
is given by $x \mapsto [\mu \mapsto x\delta_{\mu,\lambda}]$
where $\delta$ is the Kronecker delta.
Given a family of homomorphisms
$f_{\lambda}:M_{\lambda} \to N$,
define $f:M_{\infty} \to N$ by
$(x_{\lambda})_{\lambda} \mapsto \sum f_{\lambda}(x_{\lambda})$.
We see that this is the unique map.
\end{enumerate}
\end{proof}

\begin{Def}
For a complete semiring $R$,
set $R^{\prime}$ be the pre-complete
semiring, obtained from $R$ by
forgetting the infinite sum.
\begin{enumerate}
\item A \textit{\rom{(}pre-complete\rom{)} $R^{\prime}$-module} 
is a pre-complete monoid $M$
equipped with a scalar product $R^{\prime} \times M \to M$
of $R^{\prime}$, satisfying:
\begin{enumerate}[(a)]
\item $1\cdot x=x$, $0 \cdot x=0$, $a \cdot 0=0$
for any $x \in M$, $a \in R$.
\item The finite distribution law holds:
\begin{enumerate}[(i)]
\item $a(x+y)=ax+ay$ for any $a \in R$, $x,y \in M$.
\item $(a+b)x=ax+bx$ for any $a,b \in R$, $x \in M$.
\end{enumerate}
\item $(ab)x=a(bx)$ for any $a,b\in R$, $x \in M$.
\end{enumerate}
\item A map $f:M \to N$ between two $R^{\prime}$-modules
is a \textit{homomorphism of $R^{\prime}$-modules},
if it is a homomorphism of pre-complete monoids,
and preserves the scalar operation.
\item We denote by $\cat{$R^{\prime}$-pmod}$
the category of pre-complete $R^{\prime}$-modules.
\end{enumerate}
\end{Def}
\begin{Prop}
The underlying functor
$U: \cat{$R$-mod} \to \cat{$R^{\prime}$-pmod}$
has a left adjoint.
\end{Prop}
\begin{proof}
The proof is analogous
to that of Proposition ~\ref{prop:adj:completion}.
Let $M$ be a pre-complete $R^{\prime}$-module.
A \textit{$R$-filter} on $M$ is a non-empty subset $F$ of $M$
satisfying:
\begin{enumerate}
\item $x,y \in F \Leftrightarrow x+y \in F$.
\item If $a_{\lambda}x \in F$ for any $\lambda$,
then $(\sum a_{\lambda})x \in F$.
\end{enumerate}
Let $M^{\ddagger}$ be the set of 
$R$-filters on $M$.
Then, $M^{\ddagger}$ becomes a $R$-module.
Hence, we can define 
a functor $\comp^{R}:\cat{$R^{\prime}$-pmod} 
\to \cat{$R$-mod}$ by $M \mapsto M^{\ddagger}$.
The rest is the repetition of \ref{prop:adj:completion}.
\end{proof}

\begin{Rmk}
Note that $\comp^{R}(M)$
is \textit{not} algebraic, since
a infinite sum appears in the definition
of the $R$-filter.
However, the situation is different
when $R$ is algebraic.
\end{Rmk}

\begin{Def}
Let $R$ be an algebraic complete semiring.
\begin{enumerate}
\item A $R$-module $M$ is \textit{algebraic} if:
\begin{enumerate}
\item $M$ is algebraic as a complete idempotent monoid.
\item If $a \in R$, $x \in M$ are both compact,
then so is $ax$.
\end{enumerate}
\item A  homomorphism $f:M \to N$ of $R$-modules
is \textit{algebraic} if it preserves compactness.
\item We denote by $\cat{alg.$R$-mod}$
the category of algebraic $R$-modules.
\end{enumerate}
\end{Def}
If $R$ is an algebraic complete semiring
and $M$ is an algebraic $R$-module,
then the pre-complete monoid
$M_{\cpt}$ consisting of compact elements
of $M$ becomes a $R_{\cpt}$-module.
Hence we have a functor
$U_{\cpt}:\cat{alg.$R$-mod} \to \cat{$R_{\cpt}$-pmod}$.

\begin{Prop}
The above functor $U_{\cpt}$ gives an equivalence of 
categories
$U_{\cpt}:\cat{alg.$R$-mod}\simeq \cat{$R_{\cpt}$-pmod}$.
\end{Prop}
\begin{proof}
This is the analogue of Proposition
~\ref{prop:comp:isom:cpt}.
Let $M$ be a pre-complete $R_{\cpt}$-module
and $M^{\ddagger}$ be the set of
filters (\textit{not} the set of $R$-filters) of $M$.
The scalar operation on $M^{\ddagger}$
is given by
\[
F \cdot \scr{X}=\sum_{a\in F,x\in \scr{X}} \langle ax\rangle,
\]
where $F$ (resp. $\scr{X}$) is a filter
on $R_{\cpt}$ (resp. $M$).
Hence, we obtain a functor 
$\comp^{\prime}:\cat{$R_{\cpt}$-pmod} \to \cat{alg.$R$-mod}$.
It is easy to see that this functor gives the
inverse of $U_{\cpt}$.
\end{proof}

\begin{Def}
Let $R$ be a complete semiring,
and $R^{\prime}$ be the pre-complete
semiring obtained from $R$ by
forgetting the infinite sum.
Let $\sigma$ be an algebraic type
stronger than $R^{\prime}$-module.
We denote by $\sigma^{\ddagger}$
the complete algebraic type induced by $\sigma$,
extending the finite scalar operation
to the infinite scalar operation. 
\end{Def}
\begin{Prop}
Let $R$ be a complete semiring,
and $R^{\prime}$ be the pre-complete
semiring obtained from $R$ by
forgetting the infinite sum.
Let $\sigma, \tau$ be an algebraic type
stronger than $R^{\prime}$-module,
and suppose $\tau$ is stronger than $\sigma$.
Then, there is a left adjoint functor
of the underlying functor
$U_{v}:\cat{$\tau^{\ddagger}$-alg}
 \to \cat{$\sigma^{\ddagger}$-alg}$.
\end{Prop}
\begin{proof}
This is the analogue of Proposition
~\ref{prop:adj:comp:strong}.
Let
$F_{\sigma}:\cat{$\sigma $-alg} \leftrightarrows  
\cat{$\sigma^{\ddagger}$-alg}:U_{\sigma}$,
$F_{\tau}:\cat{$\tau $-alg} \leftrightarrows  
\cat{$\tau^{\ddagger}$-alg}:U_{\tau}$,
and 
$F_{\sigma}:\cat{$\sigma $-alg} \leftrightarrows  
\cat{$\tau $-alg}:U_{\sigma}$
be adjoints, respectively:
\[
\xymatrix{
\cat{$\sigma $-alg} \ar@{<->}[r]^{u} \ar@{<->}[d]^{\sigma} & 
\cat{$\tau $-alg} \ar@{<->}[d]^{\tau} \\
\cat{$\sigma^{\ddagger}$-alg} \ar@{<.>}[r]^{v} &
\cat{$\tau^{\ddagger}$-alg}
}
\]
Given a $\sigma^{\ddagger}$-algebra $A$,
set $A^{\prime}=F_{u}U_{\sigma}A \in \cat{$\tau $-alg}$.
Let $F_{v}(A)$ be the set of 
$R$-filters on $A^{\prime}$ which is infinite with respect to $A$.
Then, $F_{v}(A)$ becomes a $\tau^{\ddagger}$-algebra,
and we obtain the functor $F_{v}: \cat{$\sigma^{\ddagger}$-alg}
\to \cat{$\tau^{\ddagger}$-alg}$.
This becomes the left adjoint of $U_{v}$.
\end{proof}

\begin{Cor}
Let $R$ be a complete semiring.
The underlying functor
$\cat{$R$-alg} \to \cat{$R$-mod}$
has a left adjoint $S$.
For any $R$-module $M$, $S(M)$
is the \textit{symmetric algebra} generated by $M$. 
\end{Cor}

\subsection{Tensor products}
Throughout this subsection, we fix a complete semiring $R$.
\begin{Prop}
For any $R$-module $N$, the functor
$\Hom_{\cat{$R$-mod}}(N, -): \cat{$R$-mod} \to \cat{$R$-mod}$
has a left adjoint.
\end{Prop}
\begin{proof}
Let $M$ be another $R$-module.
A non-empty subset $F$ of $M \times N$
is a \textit{filter} if it satisfies the followings:
\begin{enumerate}
\item Let $x$, $y$ and $r$ be elements of $M$,
$N$ and $R$, respectively. Then $(rx,y) \in F$ if and only if $(x,ry) \in F$.
\item Let $x_{\lambda}$ be elements of $M$,
and $y$ be an element of $N$.
Then, $(x_{\lambda},y) \in F$ for any $\lambda$ if and only if
$(\sum_{\lambda}x_{\lambda},y) \in F$.
\item Let $y_{\lambda}$ be elements of $N$,
and $x$ be an element of $M$.
Then, $(x,y_{\lambda}) \in F$ for any $\lambda$ if and only if
$(x,\sum_{\lambda}y_{\lambda}) \in F$.
\end{enumerate}
We will denote by $\sum_{\lambda}x_{\lambda}\otimes y_{\lambda}$
the filter generated by $\{(x_{\lambda},y_{\lambda})\}_{\lambda}$.
We define $M \otimes_{R} N$ as the set of
all filters.
The supremum of $\{F_{\lambda}\}$
is the filter generated by $F_{\lambda}$'s.
For any scalar $r \in R$ and any filter $F$,
the scalar operation is defined by $r\cdot F=\{(rx,y)\}_{(x,y) \in F}$.
Given a homomorphism $f:M \to M^{\prime}$
of $R$-modules, the homomorphism
$f \otimes N:M \otimes N \to M^{\prime} \otimes N$
is defined by $\sum_{\lambda} x_{\lambda} \otimes y_{\lambda}
\mapsto \sum_{\lambda} f(x_{\lambda}) \otimes y_{\lambda}$.

Thus, the functor $- \otimes N :\cat{$R$-mod} \to \cat{$R$-mod}$
is well defined.
We will prove that this is the left adjoint of
$\Hom_{\cat{$R$-mod}}(N,-)$.
The unit $\epsilon: \Id \Rightarrow \Hom(N, - \otimes_{R} N)$
is given by $x \mapsto [y \mapsto x \otimes y]$.
The counit $\Hom(N, -) \otimes_{R} N \Rightarrow \Id$
is defined by $f \otimes y \mapsto f(y)$.
\end{proof}

\begin{Prop}
Let $R$ be a semiring,
and $A$ be a $R$-algebra.
Then, the underlying functor
$U: \cat{$R$-mod} \to \cat{$A$-mod}$
has a left adjoint $A \otimes_{R} -$.
\end{Prop}
\begin{proof}
The proof is the analogy
of classical algebras.
Let $M$ be a $R$-module.
The scalar multiplication of $A$  on $A \otimes_{R} M$
is given by $a(\sum b_{\lambda} \otimes x_{\lambda})
=\sum ab_{\lambda} \otimes x_{\lambda}$.
This gives the structure of $A$-module on $A \otimes_{R} M$.

The unit $\epsilon: \Id_{\cat{$R$-mod}} \Rightarrow U\circ (A\otimes_{R} -)$
is given by $x \mapsto 1 \otimes x$.
The counit $\eta: (A \otimes_{R}-) \circ U \Rightarrow \Id_{\cat{$A$-mod}}$
is given by
$\sum a_{\lambda} \otimes y_{\lambda} \mapsto \sum
 a_{\lambda}y_{\lambda}$.
\end{proof}
\begin{Def}
Let $F$ be the left adjoint
of the underlying functor
$U: \cat{Set} \to \cat{$R$-mod}$.
For any set $S$, $R^{S}=F(S)$
is the \textit{free $R$-module generated by $S$}.
\end{Def}

\section{Congruence relations}
\subsection{Idealic semirings}
\begin{Def}
\begin{enumerate}
\item A semiring $R$ is \textit{idealic},
if the multiplicative unit $1$ is the maximal element.
\item We denote by $\cat{IRng$^{\dagger}$}$
the full subcategory of $\cat{SRng$^{\dagger}$}$
consisting of complete idealic semirings.
\end{enumerate}
\end{Def}

\begin{Def}
Let $f:A \to B$ a homomorphism
of complete idealic semirings.
For an element $b$ of $B$,
$\pi^{-1}(b)=\sup\{x\in A \mid \pi(x) \leq b\}$
is the \textit{inverse image} of $b$.
Note that $\pi^{-1}$ neither
preserves supremums nor multipications.
\end{Def}

\subsection{Congruence relations}

\begin{Def}
Let $R$ be a complete semiring.
A \textit{congruence relation} $\mathfrak{a}$ of $R$
is an equivalence relation on $R$ satisfying the
following conditions:
\begin{enumerate}
\item If $(a_{\lambda}, b_{\lambda}) \in \mathfrak{a}$,
then $(\sum a_{\lambda}, \sum b_{\lambda}) \in \mathfrak{a}$.
\item If $(a_{i},b_{i}) \in \mathfrak{a}$ for $i=1,2$,
then $(a_{1}a_{2},b_{1}b_{2}) \in \mathfrak{a}$.
\end{enumerate}
If $\mathfrak{a}$ is a congruence relation of $R$,
then $R/\mathfrak{a}$ has a natural
structure of a complete semiring,
and there is a surjective homomorphism
$\pi: R \to R/\mathfrak{a}$ of complete semirings.
We denote by $\tilde{R}$
the set of congruence relations of $R$.
\end{Def}

The next proposition is obvious.
\begin{Prop}
Let $R$ be a complete semiring.
Then $\tilde{R}$ parametrizes
surjective homomorphisms
of complete semirings from $R$, i.e.
for any surjective homomorphism $f:R \to A$
of complete semirings,
there exists a unique congruence relation $\mathfrak{a}$
of $A$ and a natural isomorphism
$R/\mathfrak{a} \simeq A$ making
the following diagram commutative:
\[
\xymatrix{
R \ar[r]^{f} \ar[d] & A \\
R/\mathfrak{a} \ar[ru]_{\simeq}
}
\]
\end{Prop}

\begin{Prop}
The set $\tilde{R}$ has a natural
structure of an idealic $R$-algebra.
\end{Prop}
\begin{proof}
Let $\mathfrak{a}$ and $\mathfrak{b}$
be two congruence relations of $R$.
We define the multiplication $\mathfrak{a}\cdot \mathfrak{b}$
as the congruence relation generated by
$\{ (ab+a^{\prime}b^{\prime},ab^{\prime}+a^{\prime}b)\}$,
where $(a,a^{\prime}) \in \mathfrak{a}$
and $(b,b^{\prime}) \in \mathfrak{b}$.
The supremum of congruence relations
is the congruence relation generated by those.
We can easily verify that this gives
the structure of idealic semiring on $\tilde{R}$:
The multiplicative unit is $1_{\tilde{R}}=R \times R$,
and the additive unit is $0_{\tilde{R}}=\Delta$,
the diagonal subset of $R \times R$.
Finally, the homomorphism
$R \to \tilde{R}$ of semirings is given by
$a \mapsto \langle (a,0)\rangle$,
where $\langle (a,0) \rangle$ is the congruence relation
generated by $(a,0)$.
\end{proof}

\begin{Rmk}
We must notice that a
localization of a idealic semiring
(to be appeared in 4.2)
is also surjective, hence in the category of idealic semirings
we cannot distinguish localizations
from quotient rings, using only congruence relations.
\end{Rmk}

\begin{Def}
Let $R$ be a complete semiring,
$\mathfrak{a}$ be a congruence relation on $R$,
and $\pi:R \to R/\mathfrak{a}$ be
a natural map.
For an element $a$ of $R$,
$\overline{a}=\pi^{-1}\pi(a)$
is the \textit{$\mathfrak{a}$-saturation}
of $a$.
$a$ is \textit{$\mathfrak{a}$-saturated} if
$\overline{a}=a$.
\end{Def}

\begin{Lem}
\label{lem:saturation}
Let $R$ be a complete idealic semiring,
and $\mathfrak{a}$ be a congruence relation on $R$.
\begin{enumerate}
\item For any $a_{\lambda} \in R$,
$\overline{\sum \overline{a_{\lambda}}}
=\overline{\sum a_{\lambda}}$.
\item For any $a,b \in R$,
$\overline{\overline{a}\cdot \overline{b}}
=\overline{ab}$.
\end{enumerate}
\end{Lem}

\begin{proof}
\begin{enumerate}
\item
Note that $\pi(\pi^{-1}(x))=x$
for any $x \in R/\mathfrak{a}$,
since $\pi$ is surjective.
\begin{multline*}
\overline{\textstyle\sum \overline{a_{\lambda}}}
=\pi^{-1}(\pi (\textstyle\sum \overline{a_{\lambda}}))
=\pi^{-1}(\textstyle\sum \pi (\overline{a_{\lambda}})) \\
=\pi^{-1}(\textstyle\sum \pi (a_{\lambda}))
=\pi^{-1}(\pi \textstyle\sum( a_{\lambda}))
=\overline{\textstyle\sum a_{\lambda}}
\end{multline*}
\item Similar argument.
\end{enumerate}
\end{proof}

\begin{Lem}
\label{lem:semiideal:alg}
Let $R$ be an algebraic complete semiring,
and $\mathfrak{a}$ be a congruence relation of $\mathfrak{a}$.
Then, the following are equivalent:
\begin{enumerate}
\item
$R/\mathfrak{a}$ is algebraic,
and the natural map $\pi:R \to R/\mathfrak{a}$
is algebraic (When this happens, we call $\mathfrak{a}$ \textit{algebraic}.)
\item Let $a \in R$ be compact, and $b_{\lambda} \in R$.
Then, $(a+\sum b_{\lambda},\sum b_{\lambda}) \in \mathfrak{a}$
implies 
$(a+\sum^{<{\infty}} b_{\lambda},\sum^{<{\infty}} b_{\lambda})
\in \mathfrak{a}$.
\end{enumerate}
\end{Lem}
The proof is straightforward.

\begin{Def}
Let $R$ be a complete semiring.
A semiorder $\prec$ on $R$ is a 
\textit{complete \rom{(}resp. finite\rom{)} idealic semiorder} if:
\begin{enumerate}
\item $a_{\lambda} \prec b_{\lambda}$ for any $\lambda$ $\Rightarrow 
\sum a_{\lambda} \prec \sum b_{\lambda}$
(resp. $a_{1}+a_{2} \prec b_{1} +b_{2}$).
\item $a\leq b \Rightarrow a \prec b$.
\item $a_{i} \prec b_{i}$ ($i=1,2$)
$\Rightarrow a_{1}a_{2} \prec b_{1}b_{2}$.
\end{enumerate}
If $\prec$ is a complete idealic semiorder,
then the equivalence relation $\mathfrak{a}$
defined by $(a,b) \in \mathfrak{a} \Leftrightarrow a \prec b, b\prec a$ 
becomes a congruence relation on $R$.
\end{Def}

\begin{Prop}
\label{prop:semiorder:alg}
Let $R$ be an algebraic complete semiring.
Let $\prec^{f}$ be a finite idealic semiorder
satisfying:
(*) if $x$ is compact and $x \prec^{f} b$,
then there is a compact $b^{\prime} \leq b$
such that $x \prec^{f} b^{\prime}$.

Let $\prec$ be the complete idealic semiorder
generated by $\prec^{f}$.
Then, the following are equivalent:
\begin{enumerate}[(i)]
\item $a \prec b$
\item $x \prec^{f} b$ holds for any compact $x \leq a$.
\end{enumerate}
Furthermore,
if $\mathfrak{a}$ is a
congruence relation generated by $\prec$,
then $\mathfrak{a}$ is algebraic,
and $a \leq b$ in $R/\mathfrak{a}$
if and only if $a \prec b$.
\end{Prop}
\begin{proof}
(i) $\Rightarrow$ (ii) is clear.
We will show (ii) $\Rightarrow$ (i).
Let $\ll$ be the relation defined by $a \ll b$
if the condition (ii) holds.
It suffices to show that $\ll$ is actually
a complete semiorder.
Firstly, $\ll$ is a semiorder:
indeed, if $a \ll b$ and $b \ll c$,
then $x \prec^{f} b$ holds for any compact $x \leq a$.
The given condition (*) implies that
there is a compact $b^{\prime} \leq b$
such that $x \prec b^{\prime}$.
$b \ll c$ implies that $b^{\prime} \prec^{f} c$,
hence the result follows.
It is easy to see that $\ll$ is finite idealic.
It remains to show that $\ll$ is complete.
Suppose $a_{\lambda} \ll b_{\lambda}$.
Then, for any compact $x \leq \sum a_{\lambda}$,
\[
x \leq \sum^{<\infty} a_{\lambda}
\prec^{f} \sum^{<\infty} b_{\lambda}
\prec^{f} \sum b_{\lambda},
\]
hence $\sum a_{\lambda} \ll \sum b_{\lambda}$.
The rest are straightforward.
\end{proof}

\section{Topological spaces}
In this section,
we will see that the sober spaces
are the most appropriate for
considering spectrum functors.
However, most topics
in this section has been already done decades
ago in the lattice theory.

\subsection{The spectrum functor}
In the sequel,
all semirings are complete.
\begin{Def}
\begin{enumerate}
\item We denote by $\cat{Top}$,
the category of topological spaces.
\item A topological space $X$ is
\textit{sober}, if any 
(nonempty) irreducible closed subset of $X$
has a unique generic point.
We denote by $\cat{Sob}$
the full subcategory of $\cat{Top}$
consisting of 
sober topological spaces.
\end{enumerate}
\end{Def}

\begin{Prop}
The category $\cat{Sob}$
is a full co-reflective subcategory of $\cat{Top}$,
i.e. the underlying functor
$U:\cat{Sob} \to \cat{Top}$
has a left adjoint.
\end{Prop}
\begin{proof}
This proposition is standard.

For a given topological space $X$,
Let $\sob(X)$ be the set
of irreducible closed subsets of $X$.
Closed sets of $\sob(X)$ are 
of forms $V(z)=\{ c \in \sob(X) \mid c \subset z\}$,
where $z$ is a closed subset of $X$.
Given a continuous map $f:X \to Y$
between topological spaces,
$\sob(f):\sob(X) \to \sob(Y)$
is defined by
$c \mapsto \overline{f(c)}$,
where $\overline{f(c)}$
is the closure of $f(c)$.
This gives a functor
$\sob:\cat{Top} \to \cat{Sob}$.
The unit
$\epsilon:\Id_{\cat{Top}} \Rightarrow U \circ\sob$
is given by $x \mapsto \overline{\{x\}}$.
The counit 
$\eta: \sob \circ U \Rightarrow \Id_{\cat{Sob}}$
is given by $z \mapsto \xi_{z}$,
where $\xi_{z}$ is the unique generic point
of $z$.
\end{proof}
\begin{Rmk}
The counit $\eta$ in the above proof
is actually an isomorphism:
the inverse is the unit $\epsilon$.
\end{Rmk}

\begin{Def}
\begin{enumerate}
\item
Given a sober space $X$,
the set $C(X)$ of closed sets of $X$
becomes a complete idealic semiring
with idempotent multiplication:
namely, the supremum of
closed sets are their intersections,
and the multiplication of two closed sets
is the union of the two.

\item Given a continuous map $f:Y \to X$
between two sober spaces,
we can define $C(f):C(X) \to C(Y)$
by $Z \mapsto f^{-1}(Z)$.
Thus, we obtain a contravariant functor
$C :\cat{Sob}^{\op} \to \cat{IRng$^{\dagger}$}$.
\end{enumerate}
\end{Def}

\begin{Def}
Let $R$ be an idealic semiring.
\begin{enumerate}
\item An element $p$ of $R$
is \textit{prime}, if 
the subset $\{ a \in R \mid  a \nleq p\}$
is a multiplicative monoid.
\item The \textit{spectrum} $\Spec R$ of $R$
is the set of all prime elements of $R$.
For any $a \in R$,
define $V(a) \subset \Spec R$ as
the subset consisting of all prime greater than $a$:
$V(a)=\sup\{p \in \Spec R \mid a \leq p\}$.
The subset of the form $V(a)$ satisfies the
axiom of closed sets, hence gives
a topology on $\Spec R$.
It is easy to see that $\Spec R$ is
a sober space.
\item Let $f:A \to B$ be a homomorphism
of idealic semirings.
We define a continuous map $\Spec f:\Spec B \to \Spec A$
by 
\[
y \mapsto f^{-1}(y)=\sup \{ x \in R \mid f(x) \leq y\}.
\]
Hence, we obtain a contravariant functor
$\Spec:\cat{IRng$^{\dagger}$} \to \cat{Sob}^{\op}$,
to which we refer as the \textit{spectrum functor}.
\end{enumerate}
\end{Def}
Note that $\Spec f$ is well defined.
Indeed, let $p \in B$ be a prime element.
Let $a, b \in A$ be two elements
satisfying $a,b \nleq f^{-1}(p)$.
This means that $f(a),f(b) \nleq p$,
hence $f(ab) \nleq p$ since $p$ is prime.
Thus we obtained $ab \nleq f^{-1}(p)$.

\begin{Prop}
\label{prop:top:spec:adjoint}
The spectrum functor $\Spec$ is the left adjoint of
$C:\cat{Sob}^{\op} \to \cat{IRng$^{\dagger}$}$.
\end{Prop}
\begin{proof}
The unit $\epsilon:\Id_{\cat{IRng$^{\dagger}$}} \Rightarrow C \circ \Spec$
is given by $A \ni a \mapsto V(a)$.
The counit $\eta: \Spec \circ C \Rightarrow \Id_{\cat{Sob}^{\op}}$
is given by 
\[
X \ni x \mapsto \overline{\{x \}} \in \Spec \circ C(X).
\]
\end{proof}
\begin{Rmk}
The above $\eta$ is actually a homeomorphism:
the inverse is given by sending
an irreducible closed set $Z$ to its unique
generic point.
This shows that
$\cat{Sob}^{\op}$ is a reflective full
subcategory of $\cat{IRng$^{\dagger}$}$.
\end{Rmk}

\begin{Def}
\begin{enumerate}
\item
Let $X$ be a sober space.
We call $X$ \textit{algebraic} if:
\begin{enumerate}
\item $X$ is quasi-compact.
\item $X$ is quasi-seperated,
i.e. for any two quasi-compact open subsets
$U,V$ of $X$, $U \cap V$ is quasi-compact
(cf. ~\cite{EGA},Chap I. 1.2).
\item Any open subset of $X$ is a
union of some quasi-compact open subsets.
\end{enumerate}
\item A \textit{morphism $f:X \to Y$
of algebraic sober spaces}
is a continuous map which is \textit{quasi-compact}:
the inverse image of a quasi-compact open subset
is again quasi-compact.
\item We denote by $\cat{alg.Sob}$
the category of algebraic sober
spaces.
\end{enumerate}
\end{Def}
A sober space $X$ is algebraic if and only if
$C(X)$ is algebraic.
An algebraic sober space is a quasi-compact
coherent space, and vice versa.

\begin{Prop}
\label{prop:spec:alg}
Let $R$ be an algebraic complete idealic semiring,
and $X=\Spec R$. Denote by $D(a)$
the open subset $X\setminus V(a)$
for any element $a \in R$.
\begin{enumerate}
\item Let $U$ be
an open subset of $X$.
Then, $U$ is quasi-compact if and only if 
there is a compact element $a$ of $R$
such that $U=D(a)$.
\item $X$ is algebraic.
\item If $\varphi: B \to A$
is a homomorphism of algebraic
complete idealic semirings,
then $\Spec(\varphi):\Spec B \to \Spec A$
is algebraic.
\end{enumerate}
\end{Prop}
\begin{proof}
\begin{enumerate}
\item Suppose $U$ is quasi-compact.
There is an element $b \in R$
such that $U=D(b)$.
There is a covering $b=\sum b_{\lambda}$
of $b$ by compact elements.
Then $U=\cup D(b_{\lambda})$,
and since $U$ is quasi-compact,
there is a finite subcover $U=\cup^{<\infty} D(b_{\lambda})$.
This means that $U=D(\sum^{<\infty}b_{\lambda})$,
and $\sum^{<\infty}b_{\lambda}$ is compact.
The converse is easy.
\item
First, quasi-compactness of $X$ follows from (1)
since the unit $1$ is compact.
Let $U_{1},U_{2}$ be two quasi-compact
open subsets of $X$.
Then, there are compact elements $a_{i} \in R$
such that $U_{i}=D(a_{i})$.
Then $a_{1}a_{2}$ is compact,
so $U_{1} \cap U_{2}=D(a_{1}a_{2})$
is quasi-compact. Thus,
$X$ is quasi-seperated.
Any open subset of $X$ 
is a union of quasi-open subsets,
since any element of $R$ is algebraic.
Thus, $X$ is algebraic.
\item Easy.
\end{enumerate}
\end{proof}

\begin{Cor}
\label{Cor:adj:alg:irng:sob}
The spectrum functor
gives a functor 
$\Spec: \cat{alg.IRng$^{\dagger}$} 
\to \cat{alg.Sob}^{\op}$
when restricted to $\cat{alg.IRng$^{\dagger}$}$,
and is the left adjoint of
$C:\cat{alg.Sob}^{\op} \to \cat{alg.IRng$^{\dagger}$}$.
\end{Cor}
\begin{proof}
To see the adjointness,
it suffices to show that the
unit $\epsilon$ and the counit $\eta$
in Proposition ~\ref{prop:top:spec:adjoint}
are algebraic.
The above proposition shows that
$\epsilon$ is algebraic, and $\eta$
is algebraic since it is a homeomorphism.
\end{proof}
We will strengthen this
result later.
(See Proposition ~\ref{prop:equiv:sob:irng}.)

\subsection{Localization}

In this subsection,
we gather some basic facts
on localization of idealic semirings.
Most of the results are similar to
those of ordinary rings.
\begin{Def}
Let $R$ be a complete idealic semiring.
\begin{enumerate}
\item A subset $\Sigma$ of $R$
is a \textit{multiplicative system of $R$}
if it is a submonoid of $R$.
$\Sigma$ is \textit{compact},
if it consists of compact elements.
\item The \textit{localization of $R$ along $\Sigma$}
is the homomorphism $R \to R/\mathfrak{a}$
of complete idealic semirings, where $\mathfrak{a}$
is a congruence relation on $R$ generated by pairs $(1,s)$,
where $s \in \Sigma$.
We write $\Sigma^{-1}R$ instead
of $R/\mathfrak{a}$.
\item If $\Sigma_{f}=\{f^{n}\}_{n \in \NNN}$
for an element $f \in R$,
we set $R_{f}=\Sigma_{f}^{-1}R$.
\item If $\Sigma_{p}=\{ x \mid x \nleq p\}$
for a prime element $p$ of $R$,
we set $R_{p}=\Sigma_{p}^{-1}R$,
and refer to this semiring
the \textit{localization of $R$ along $p$}.
If $R$ is algebraic, then
we set $\Sigma_{p}=\{x \mid \rom{$x$ is compact and $x \nleq p$}\}$
instead.
\item Given an element $a$ of $R$,
$\overline{a}=\pi^{-1}\pi(a)$
is the \textit{$\Sigma$-saturation of $a$},
where $\pi:R \to \Sigma^{-1}R$
is the natural map.
$a$ is \textit{$\Sigma$-saturated} if
$\overline{a}=a$.
\end{enumerate}
\end{Def}

\begin{Lem}
\label{lem:loc:relation}
Let $R$ be an algebraic complete idealic semiring,
and $\Sigma$ be a compact multiplicative system.
Let $\mathfrak{a}$ be the congruence relation
defined by $\Sigma$.
Then $(f,g) \in \mathfrak{a}$
if and only if:

($\mathfrak{b}$) For any compact $x$ smaller than $f$,
there exists $s\in \Sigma$ such that
$sx \leq g$.
For any compact $y$ smaller than $g$,
there exists $t\in \Sigma$ such that
$ty \leq f$.
\end{Lem}
\begin{proof}
We make use of Proposition ~\ref{prop:semiorder:alg}.
Let $\prec^{f}$ be a relation defined by
$a \prec^{f}b$ if:

For any compact $x \leq a$,
there exists $s \in \Sigma$ such that $sx \leq b$.

We can easily see that $\prec^{f}$
is a finite idealic semiorder, and $\mathfrak{a}$
is the congruence relation generated by $\prec^{f}$.
It remains to show that
if $x \prec^{f} b$ is compact,
then there exists a compact $b^{\prime} \leq b$
such that $x \prec^{f} b^{\prime}$.
Since $sx \leq  b$ for some $s \in \Sigma$,
and $sx$ is compact,
there exists a compact $b^{\prime} \leq b$
satisfying $sx \leq b^{\prime}$.
Hence, the result follows.
\end{proof}

\begin{Cor}
\label{cor:alg:cpt:loc:alg}
Let $R$ be an algebraic
complete idealic semiring,
and $\Sigma$ be a compact multiplicative
system. Then:
\begin{enumerate}
\item $\Sigma^{-1}R$ is algebraic
and the natural map 
$\pi: R \to \Sigma^{-1}R$ is also algebraic.
\item More generally, let $\Sigma_{1}$, $\Sigma_{2}$
be compact multiplicative systems.
If there is a homomorphism
$f:\Sigma_{1}^{-1}R \to \Sigma_{2}^{-1}R$
of $R$-algberas,
then $f$ is algebraic.
\item For any $a \in R$,
$\overline{a}$ is the supremum of all compact $x \in R$
satisfying $sx \leq a$ for some $s \in \Sigma$.
\end{enumerate}
\end{Cor}
\begin{proof}
These are direct consequences
of Lemma ~\ref{lem:saturation} and Lemma ~\ref{lem:loc:relation}.
Let us prove (2).
Suppose $x \in \Sigma^{-1}R$
is a compact element.
Then, we can find a compact $a \in R$
such that $\pi_{1}(a)=x$,
where $\pi_{i}:R \to \Sigma^{-1}_{i}R$
is the natural map:
\[
\xymatrix{
R \ar[d]_{\pi_{1}} \ar[rd]^{\pi_{2}} \\
\Sigma_{1}^{-1}R \ar[r]_{f} & \Sigma_{2}^{-1}R
}
\]
Then, $f(x)=\pi_{2}(a)$, hence compact.

\end{proof}

\begin{Cor}
\label{Cor:top:property:spec}
Let $R$ be an algebraic complete
idealic semiring.
\begin{enumerate}
\item Let $f$ be a compact element in $R$.
The localization map $R \to R_{f}$
induces $\Spec R_{f} \to \Spec R$,
which is an open immersion:
$\Spec R_{f} \simeq D(f)=\Spec R\setminus V(f)$.
\item The open sets of the form $D(f)$
with $f$ compact
gives an open basis of $\Spec R$.
\end{enumerate}
\end{Cor}
\begin{Cor}
Let $R$ be an algebraic complete idealic semiring,
and $p$ a prime element of $R$.
Then, $R_{p}$ is an algebraic complete
idealic local semiring, with $pR_{p}$
being the unique maximal ideal.
Here, $pR_{p}$ is the image of
$p$ via the natural map $\pi:R \to R_{p}$.
\end{Cor}
\begin{Lem}
\label{lem:alg:exist:prime}
Let $R$ be a complete idealic semiring.
\begin{enumerate}
\item Any maximal element (of $R\setminus \{1\}$)
is prime.
\item Suppose $R$ is algebraic.
Then, for any non-unit element $a \neq 1$,
there exists a maximal element larger than $a$.
\end{enumerate}
\end{Lem}
\begin{proof}
\begin{enumerate}
\item
Let $m$ be a maximal element of $R$,
and suppose $a \nleq m$ and $b \nleq m$.
then $a+m=b+m=1$ since $m$ is maximal.
Hence, $m+ab \geq (m+a)(m+b)=1$
which shows that $ab \nleq m$.
Thus, $m$ is prime.
\item Let 
$\scr{S}$ be a set of non-unit elements
which are larger than $a$.
Then $\scr{S}$ is a inductively ordered
set, since $1$ is compact.
Hence, there is a maximal element of $\scr{S}$
by Zorn's lemma.
\end{enumerate}
\end{proof}

\begin{Lem}
\label{lem:alg:v:equiv}
Let $R$ be an algebraic complete idealic ring,
and $a,b \in R$.
Then, the following are equivalent:
\begin{enumerate}[(i)]
\item For any compact $x \leq a$,
there exists a natural number $n$ such that $x^{n} \leq b$.
\item $\sqrt{a} \leq \sqrt{b}$,
where $\sqrt{a}=
\sup\{ x \mid \text{$x^{n} \leq a$ for some $n$}\}$.
\item $V(a) \supset V(b)$ on $\Spec R$.
\end{enumerate}
Further, if $a$ is compact,
these are also equivalent to:
\begin{enumerate}[(i)]
\setcounter{enumi}{3}
\item $\mathfrak{a} \geq \mathfrak{b}$,
where $\mathfrak{a}$ (resp. $\mathfrak{b}$)
is a congruence relation generated by $(1,a)$ (resp. $(1,b)$).
\end{enumerate}
\end{Lem}
\begin{proof}
(i)$\Rightarrow$ (ii):
If $x \leq \sqrt{a}$ is any compact element,
then there exists a natural number $m$
such that $x^{m} \leq a$.
Then, (i) implies that $x^{mn} \leq b$
for some $n$, hence $x \leq \sqrt{b}$.
This shows that $\sqrt{a}\leq \sqrt{b}$.

(ii)$\Rightarrow$ (iii):
If $p \in V(b)$, then $p \geq b$.
Let $x \leq a$ be any compact element.
Then $x \leq \sqrt{a} \leq \sqrt{b}$ implies
$x^{m} \leq b$ for some $m$.
Since $x^{m} \leq p$ and $p$ is prime, we have
$x \leq p$. Therefore, $a \leq p$ and $p \in V(a)$.

(iii)$\Rightarrow$ (i):
Suppose there exists a compact $x \leq a$
such that $x^{n} \nleq b$ for any $n$.
Set $A=R/b$: this is the quotient ring of $R$
divided by a congruence relation generated by $(b,0)$.
Then $x^{n} \nleq 0$
for all $n$ in $A$.
This shows that $A_{x} \neq 0$.
Since $A_{x}$ is algebraic, we can find a prime
element $p$ of $A_{x}$, by Lemma ~\ref{lem:alg:exist:prime}.
Let $\varphi:R \to A_{x}$ be the canonical map.
Then, $\varphi^{-1}(p)$ is contained in $V(b)$,
but $\varphi^{-1}(p) \notin V(a)$,
since $\varphi^{-1}(p) \ngeq x$.
This is a contradiction.

Suppose $a$ is compact.

(i) $\Rightarrow$ (iv):
$a^{n} \leq b$ holds for some $n$,
and since $(a^{n},1) \in \mathfrak{a}$,
we have $(b,1) \in \mathfrak{a}$.

(iv) $\Rightarrow$ (i):
$(1,b) \in \mathfrak{b}$,
so $(1,b) \in \mathfrak{a}$ by assumption.
This is equivalent to $a^{n} \leq b$ for some $n$.
\end{proof}

\begin{Cor}
Let $R$ be an algebraic complete idealic semiring,
and $U$ be a quasi-compact open
subset of $X=\Spec R$. Set $Z=X\setminus U$.
Let $f,g \in R$ be compact elements satisfying
$V(f)=V(g)=Z$.
Then $R_{f} \simeq R_{g}$.
\end{Cor}

\begin{Def}
Let $R$ be an algebraic complete
idealic semiring.
By the above corollary,
we may define a \textit{localization
$R_{Z}$ of $R$ along any compact $Z \in C(X)$},
by $R_{Z}=R_{f}$ for any compact $f \in R$
such that $V(f)=Z$.
Also, for a non-compact $Z \in C(X)$,
set $R_{[Z]}=\underleftarrow{\lim}_{Z^{\prime} \leq Z}R_{Z^{\prime}}$,
where $Z^{\prime}$ runs through all the
compact elements smaller than $Z$, and the limit is taken
within the category of \textit{algebraic} semirings.
\end{Def}

\begin{Rmk}
Note that $R_{[Z]}$ is not isomorphic
to $R_{Z}$ in general: if $Z$ is not compact,
$R_{Z}$ may not be even algebraic.
\end{Rmk}

The next is the key lemma,
which is indispensable when
constructing idealic schemes.
\begin{Lem}
\label{lem:patch:cond:idealic}
Let $R$ be an algebraic complete
idealic semiring, and 
\linebreak $s,s_{1},\cdots, s_{n}$
be compact elements satisfying $s=\sum s_{i}$.
If there are elements $f_{i} \in R_{i}=R_{s_{i}}$
satisfying $f_{i}=f_{j}$ in $R_{ij}=R_{s_{i}s_{j}}$,
then there is a unique element $f \in R_{s}$
such that $f=f_{i}$ in $R_{s_{i}}$.
\end{Lem}
\begin{proof}
First, we will show the uniqueness of $f$.
Let $f$ and $g$ be elements of $R_{s}$,
such that $f=f_{i}=g$ in $R_{i}$.
For any compact $x \leq f$,
there exist natural numbers $m_{i}$
such that $s_{i}^{m_{i}}x \leq g$ in $R$
for any $i$.
Set $m=\max_{i} m_{i}$.
Then, $\sum s_{i}=s$ implies
\[
\sum s_{i}^{m_{i}} \geq \sum s_{i}^{m} \geq (\sum s_{i})^{nm}
=s^{nm}
\]
so that $s^{nm}x \leq g$.
This means that $x \leq g$ in $R_{s}$, so $f \leq g$.
$g \leq f$ can be shown in a similar way.

Next, we will show the existence.
Set $f=\sum f_{i}$.
To show $f=f_{i}$ in $R_{i}$,
it suffices to prove $f_{j} \leq f_{i}$ in $R_{i}$.
Since $f_{i}=f_{j}$ in $R_{ij}$,
there exists natural numbers $m_{ij}$ such that
$(s_{i}s_{j})^{m_{ij}}f_{j} \leq f_{i}$
for any $i,j$.
Set $m=\max_{i,j} m_{ij}$. Then
\[
\sum_{j} (s_{i}s_{j})^{m_{ij}}
\geq \sum_{j} (s_{i}s_{j})^{m}
\geq s_{i}^{m}(\sum_{j} s_{j})^{mn}
= s_{i}^{m}s^{mn}
\]
so that $s_{i}^{m}s^{mn}f_{j} \leq f_{i}$.
Since $s_{i} \leq s$, this means that
$s_{i}^{m(n+1)}f_{j} \leq f_{i}$,
hence $f_{j} \leq f_{i}$ in $R_{i}$.
\end{proof}
\begin{Rmk}
Note that this lemma also
holds for pre-complete idealic semirings.
\end{Rmk}

\subsection{Comparison with the Stone-\v{C}ech compactification}

\begin{Def}
We denote by $\cat{alg.IIRng$^{\dagger}$}$
the full subcategory of \linebreak
$\cat{alg.IRng$^{\dagger}$}$
consisting of algebraic complete idealic semirings
with idempotent multiplications.
\end{Def}

\begin{Prop}
\label{prop:equiv:sob:irng}
The spectrum functor
and the functor $C$ introduced in
~\ref{Cor:adj:alg:irng:sob}
gives an equivalence
between the category $\cat{alg.IIRng$^{\dagger}$}$
and \linebreak
$\cat{alg.Sob}^{\op}$, the (opposite) category
of algebraic sober spaces.
\end{Prop}
\begin{proof}
We only need to prove that the unit 
$\epsilon:\Id_{\cat{alg.IIRng$^{\dagger}$}} \Rightarrow C \circ \Spec$
of Corollary ~\ref{Cor:adj:alg:irng:sob}
is a natural isomorphism.
Let $R$ be an algebraic complete idealic
semiring, and $a,b  \in R$.
It suffices to show that $V(a) \leq V(b)$
implies $a \leq b$.
If $V(a) \leq V(b)$, then
Lemma ~\ref{lem:alg:v:equiv}
says that for any compact $x$,
$x^{n} \leq b$ for some $b$,
but since the multiplication is idempotent,
this means that $x \leq b$, hence $a \leq b$.
\end{proof}

Making use of Corollary
~\ref{cor:alg:sigma:alg:comp},
we have the following:

\begin{Prop}
The underlying functor
$U: \cat{alg.Sob}\to \cat{Sob}$
has a left adjoint $\alg$.
\end{Prop}
\begin{proof}
The functor $\alg$ preserves
the idempotency of multiplication,
$\cat{Sob}$ is a subcategory of $\cat{IIRng$^{\dagger}$}$,
and $\cat{alg.Sob}$ coincides with
$\cat{alg.IIRng$^{\dagger}$}$,
which shows that $\alg$
sends a sober space to an algebraic sober space.

This is enough for a proof,
but let us write it down this functor
explicitly, for future references.

Given a sober space $X$,
let $\alg(X)$ be the set
of prime filters on $C(X)$:
a prime filter $F$ is a non-empty subset
of $C(X)$ satisfying:
\begin{enumerate}
\item $C_{1},C_{2} \in F \Leftrightarrow C_{1}\cap C_{2} \in F$.
\item $C_{1}, C_{2} \notin F \Rightarrow C_{1} \cup C_{2} \notin F$.
\end{enumerate}
A closed set of $\alg(X)$ is of a form
$V(a)=\{ F \mid  a \subset F\}$,
where $a$ is a filter on $C(X)$.
Given a continuous map $f:X \to Y$
of sober spaces, a continuous map
$\alg(f):\alg(X) \to \alg(Y)$ is given by
\[
\alg(X) \ni p \mapsto \sum_{f^{-1}(z) \in p}\langle z \rangle \in \alg(Y),
\]
where $z$ runs through all
the closed subsets of $Y$,
satisfying $f^{-1}(z) \in p$.
Thus, we have a functor $\alg:\cat{Sob} \to \cat{alg.Sob}$.
The unit $\epsilon: \Id_{\cat{Sob}} \Rightarrow U\alg$
of the adjoint is given by $x \mapsto \langle \overline{\{x\}} \rangle$,
and the counit $\eta:\alg U \Rightarrow \Id_{\cat{alg.Sob}}$
is given by $p \mapsto \cap_{c\in p} c$.
\end{proof}

\begin{Prop}
\label{prop:haus:alg:preserve}
If $X \in \cat{Top}$ is a Hausdorff space,
then $\alg(X)$ is also Hausdorff.
\end{Prop}
\begin{proof}
The topological space $X$ is Hausdorff if
and only if the diagonal functor $\Delta:X \to X \times X$
is a closed immersion.
It is obvious that the functor $\alg$ preserves
closed immersion, hence $\alg(X) \to \alg(X) \times \alg(X)$
is also a closed immersion.
\end{proof}

\begin{Def}
Hausdorff algebraic sober spaces are called
\textit{Stone spaces}.
The category of Stone spaces and continuous maps
is denoted by $\cat{Stone}$.
\end{Def}
\begin{Rmk}
Note that a continuous map between Stone spaces
are already quasi-compact:
an open subset $U$ of a Stone space $X$
is quasi-compact if and only iff $U$ is clopen,
since $X$ is Hausdorff.
\end{Rmk}

\begin{Def}
An II-Ring $R$ is a \textit{Boolean algebra},
if there is a unary operator $\neg:R \to R$
such that
\[
a+\neg a=1, \quad  a \cdot \neg a=0.
\]
\end{Def}
\begin{Prop}
\label{prop:stone:alg:adjoint}
\begin{enumerate}
\item The category $\cat{Bool}$ of Boolean algebras
is equivalent to the opposite category of $\cat{Stone}$.
\item The underlying functor $\cat{Stone} \to \cat{alg.Sob}$
has a left adjoint and a right adjoint.
\end{enumerate}
\end{Prop}
\begin{proof}
\begin{enumerate}
\item This is obvious from the above remark.
\item The existence of the right adjoint follows from
the fact that we have a left adjoint of the underlying functor
$\cat{alg.IIRng$^{\dagger}$} \to \cat{Bool}$.

We will construct the left adjoint.
Let $X$ be an algebraic sober space,
and $A=C(X)_{\cpt}$ be the corresponding pre-complete 
idealic semiring.
Let $B$ be the subring of $A$ consisting of 
elements with \textit{negation}, namely elements $x \in A$
which satisfies $xy=0$ and $x+y=1$ for some $y \in A$.
Note that $y$ is uniquely determined by $x$.
Therefore, we see that $B$ becomes a Boolean algebra.
This correspondence gives a functor $G:\cat{alg.IIRng$^{\dagger}$}
 \to \cat{Bool}$,
and hence $G^{\op}: \cat{alg.Sob} \to \cat{Stone}$.
For any Boolean algebra $C$,
any morphism $C \to A$ of pre-complete idempotent
semirings factors through $C \to B$.
This shows that $G^{\op}$ is the left adjoint
of the underlying functor.
\end{enumerate}
\end{proof}

\begin{Rmk}
Recall the Stone-\v{C}ech compactification:
\begin{Thm}
\label{thm:stone:cech}
Let $U:\cat{CptHaus} \to \cat{Haus}$
be the underlying functor from
the category $\cat{CptHaus}$
of compact Hausdorff
spaces to the category $\cat{Haus}$
of Hausdorff spaces.
Then, $U$ has a left adjoint $\beta$.
The unit morphism $X \to \beta X$
is an open immersion if and only if $X$
is locally compact Hausdorff.
\end{Thm}
Note that the underlying functor $\cat{Stone} \to \cat{Haus}$
factors through $\cat{CptHaus}$.
This implies that, $X \to \alg(X)$ factors
through $\beta(X)$ for any Hausdorff space $X$.
Usually, $\alg(X)$ does not coincide with $\beta(X)$.

\end{Rmk}
\begin{Exam}
\begin{enumerate}
\item If $X$ is a discrete space,
then $\beta(X)$ coincides with $\alg(X)$:
indeed, $C(X)=\prod_{x \in X} \FF_{1}$,
and there is a bijection between
the set of points of $\Spec C(X)_{\cpt}$
and set of ultrafilters on $X$.
This also coincides with the set of points
of $\beta(X)$. We can also see that
$\beta(X) \to \alg(X)$ is in fact, a homeomorphism.
\item The real line $\RR$ is a locally compact Hausdorff
space, hence it becomes an open subset of $\beta\RR$.
On the other hand, $\alg(\RR)$ is a one point set,
since $\RR$ is connected, but Stone spaces
are totally disconnected.
\end{enumerate}
\end{Exam}

\section{Schemes}
\subsection{Sheaves}

When considering sheaves
of complete-algebra valued,
we need some special care,
for there are some obstructions
when applying the sheaf theory
to algbras admitting infinite operations
(in our case, the supremum map.):
we don't have sheafifications in general,
and the stalk of such sheaves do not
admit a natural induced infinite operators.
Thus, the notion of algebraicity is essential
in the following.

\begin{Def}
Let $X$ be a topological space.
\begin{enumerate}
\item We regard the preordered set $C(X)$ as 
a category.
\item Let $\scr{A}$ be a category.
A functor $\scr{F}:C(X)^{\op} \to \scr{A}$
is a \textit{$\scr{A}$-valued presheaf on $X$}.
A morphism of presheaves on $X$ is a natural transformation.
A presheaf is a \textit{sheaf} if $\scr{F}$ is a continuous functor.
\item We denote by $\cat{$\scr{A}$-PSh/$X$}$
(resp. $\cat{$\scr{A}$-Sh/$X$}$) the category
of $\scr{A}$-valued presheaves (resp. sheaves) on $X$.
\end{enumerate}
\end{Def}
Here, we chose the definition
of a sheaf to be a contravariant functor
from the category $C(X)$ of closed subsets of $X$:
this is just for convenience sake.

\begin{Def}
Let $\sigma$ be a pre-complete algebraic type
and $X$ be a topological space.
Then, the underlying functor
$U:\cat{$\cat{alg.$\sigma^{\dagger}$-alg}$-Sh/$X$}
\to \cat{$\cat{alg.$\sigma^{\dagger}$-alg}$-PSh/$X$}$
has a left adjoint $S$,
which we call the \textit{sheafification}.
\end{Def}
\begin{proof}
This is obvious,
since $\cat{alg.$\sigma^{\dagger}$-alg}$
is equivalent to $\cat{$\sigma $-alg}$
by Proposition ~\ref{prop:comp:isom:cpt},
and $\cat{$\sigma $-alg}$-valued presheaves
admit sheafifications.
\end{proof}
Next, we will see what happens
when the underlying topological space
is algebraic.
\begin{Def}
Let $X$ be an algebraic sober space.
Then, $C(X)_{\cpt}$ becomes
a pre-complete idealic semiring
with idempotent multiplication.
Let $\scr{A}$ be any category.
\begin{enumerate}
\item We call a functor $C(X)_{\cpt}^{\op} \to \scr{A}$
a \textit{$\scr{A}$-valued presheaf on $X^{\cpt}$}.
It is a \textit{sheaf} if it is finite continuous,
i.e. preserves fiber products.
\item We denote by $\cat{$\scr{A}$-PSh/$X^{\cpt}$}$
(resp. $\cat{$\scr{A}$-Sh/$X^{\cpt}$}$)
the category of $\scr{A}$-valued presheaves
(resp. sheaves) on $X^{\cpt}$.
\end{enumerate}
\end{Def}
\begin{Prop}
\label{prop:equiv:sheaf:alg}
Let $X$ be an algebraic sober space,
and $\scr{A}$ be a small complete category.
Then, the underlying functor
$U^{\cpt}:\cat{$\scr{A}$-Sh/$X$} \to \cat{$\scr{A}$-Sh/$X^{\cpt}$}$
is an equivalence of categories.
\end{Prop}
\begin{proof}
Let $\scr{F}$ be a $\scr{A}$-valued sheaf on $X^{\cpt}$.
We define a $\scr{A}$-valued sheaf $\scr{F}^{\dagger}$
on $X$ by
\[
\scr{F}^{\dagger}(Z)=\underleftarrow{\lim}_{Z^{\prime} \leq Z}
\scr{F}(Z^{\prime}),
\]
where $Z^{\prime}$ runs through all the compact
elements smaller than $Z$.
We claim that $\scr{F}^{\dagger}$ is indeed a sheaf.

For any $Z \in C(X)$, any covering $Z=\sum Z_{\lambda}$
and $a_{\lambda} \in \scr{F}^{\dagger}(Z_{\lambda})$ satisfying
$a_{\lambda}|_{Z_{\lambda}Z_{\mu}}
=a_{\mu}|_{Z_{\lambda}Z_{\mu}}$,
 we will show that there exists a unique $a \in \scr{F}^{\dagger}(Z)$
 such that $a|_{Z_{\lambda}}=a_{\lambda}$.
We have $a_{\lambda}|_{WC_{\lambda}C_{\mu}}
=a_{\mu}|_{WC_{\lambda}C_{\mu}}$ for
any compact $W \leq Z$ and $C_{\lambda} \leq Z_{\lambda}$.
Since $\scr{F}$ is a sheaf
and $W$ is covered by finitely many $WC_{\lambda}$'s,
there is a unique $a_{W} \in \scr{F}^{\dagger}(W)$ such that
$a_{W}|_{WC_{\lambda}}=a_{\lambda}|_{WC_{\lambda}}$.
By the definition of $\scr{F}^{\dagger}$, the $a_{W}$'s patch together
to give a section $a \in \scr{F}^{\dagger}(Z)$.
It is clear that $a|_{Z_{\lambda}}=a_{\lambda}$.
Uniqueness of $a$ is also clear from the construction.

Given a morphism $f:\scr{F} \to \scr{G}$
of sheaves on $X^{\cpt}$,
we define $f^{\dagger}:\scr{F}^{\dagger} \to \scr{G}^{\dagger}$ by
\[
\scr{F}^{\dagger}(Z)=
\underleftarrow{\lim}_{Z^{\prime}} \scr{F}(Z^{\prime})
\stackrel{f}{\to}
\underleftarrow{\lim}_{Z^{\prime}} \scr{G}(Z^{\prime})
=\scr{G}^{\dagger}(Z).
\]
This is well defined.
Hence, we have a functor
$\comp: \cat{$\scr{A}$-Sh/$X^{\cpt}$} \to
\cat{$\scr{A}$-Sh/$X$}$.

We will see that this is the left adjoint of $U^{\cpt}$.
The unit $\epsilon:\Id_{\cat{$\scr{A}$-Sh/$X^{\cpt}$}}
\Rightarrow  U^{\cpt}\comp$ is given by the natural isomorphism
$\scr{F}(Z) \simeq \scr{F}^{\dagger}(Z)$
for any compact $Z$.
The counit  $\eta:\comp U^{\cpt} \to 
\Id_{\cat{$\scr{A}$-Sh/$X$}}$ is given by the natural isomorphism
\[
(\scr{G}^{\cpt})^{\dagger}(Z)
=\underleftarrow{\lim}_{Z^{\prime} \leq Z:\rom{cpt}}
\scr{G}(Z^{\prime}) \simeq \scr{G}(Z).
\]
\end{proof}

\subsection{Idealic schemes}

Here, we will introduce
a notion of idealic schemes.

The significant difference between
idealic schemes and the usual scheme
is that we can construct
a universal idealic scheme from
its underlying space.

\begin{Prop}
Let $X$ be a topological space,
$Z$ be a closed subset of $X$,
and $U=X\setminus Z$ be the complement of $Z$.
Then, the restriction homomorphism
$\pi:C(X) \to C(U)$
induces an isomorphism $C(X)_{Z} \simeq C(U)$,
where $C(X)_{Z}$ is the localization along $Z$.
\end{Prop}
\begin{proof}
Since $\pi(Z)=1$ in $C(U)$,
$\pi$ factors through $C(X)_{Z}$:
\[
\xymatrix{
C(X) \ar[r]^{\pi} \ar[d] & C(U) \\
C(X)_{Z} \ar[ru]_{\tilde{\pi}}
}
\]
It is clear that $\pi$ is surjective,
hence so is $\tilde{\pi}$.
It remains to show that $\tilde{\pi}$ is injective.
Let $\mathfrak{a}$ be the congruence relation
generated by $(1,Z)$. Then we see that
\[
\mathfrak{a}=\{ (a,b) \in C(X) \times C(X) \mid aZ=bZ\}.
\]
Indeed, the righthand side contains $(1,Z)$, and
it is easy to see that this is a congruence relation,
since the multiplication is idempotent.
If $\pi(a)=\pi(b)$, then $aZ=a \cup Z=b \cup Z=bZ$,
hence $a$ coincides with $b$ in $C(X)_{Z}$.
\end{proof}

\begin{Def}
A functor 
$\tau_{X}:C(X)^{\op} \to \cat{IRng$^{\dagger}$}$
is a sheaf on $X$ defined by 
$Z \mapsto C(X)_{Z} \simeq C(X \setminus Z)$.
\end{Def}

\begin{Def}
Let $\alpha: \cat{IRng$^{\dagger}$}
 \to \cat{IIRng$^{\dagger}$}$ be the left adjoint
 of the underlying functor 
$ \cat{IIRng$^{\dagger}$} \to \cat{IRng$^{\dagger}$}$.
\begin{enumerate}
\item A triple $X=(X, \scr{O}_{X},\beta_{X})$
is an \textit{idealic scheme} if:
\begin{enumerate}
\item $X$ is a sober space.
\item $\scr{O}_{X}$ is $\cat{IRng$^{\dagger}$}$-valued sheaf.
\item $\beta_{X}:\alpha\scr{O}_{X} \to \tau_{X}$ is a morphism of 
$\cat{IIRng$^{\dagger}$}$-valued presheaves.
\item The restriction maps reflect localization:
for any $W \leq Z$
of closed subsets of $X$,
and any $a \in \scr{O}_{X}(Z)$
satisfying $\beta_{X}(\alpha(a)) \geq W$,
the restriction functor $\scr{O}_{X}(Z) \to \scr{O}_{X}(W)$
factors through the localization
$\scr{O}_{X}(Z)_{a}$. 
\[
\xymatrix{
\scr{O}_{X}(Z) \ar[r]^{\rom{res}} \ar[d] &
\scr{O}_{X}(W) \\
\scr{O}_{X}(Z)_{a} \ar[ru]
}
\]
\end{enumerate}
\item A \textit{morphism 
$(X,\scr{O}_{X},\beta_{X}) \to (Y,\scr{O}_{Y},\beta_{Y})$
of idealic schemes}
is a pair $(f,f^{\#})$ such that:
\begin{enumerate}
\item $f:X \to Y$ is a continuous map.
\item $f^{\#}:\scr{O}_{Y} \to f_{*}\scr{O}_{X}$
is a morphism of $\cat{IRng$^{\dagger}$}$-valued
sheaves on $Y$.
\item The following diagram of $\cat{IIRng$^\dagger $}$-valued
presheaves is commutative:
\[
\xymatrix{
\alpha\scr{O}_{Y} \ar[d]_{\beta_{Y}} \ar[r]^{\alpha f^{\#}} &
\alpha f_{*}\scr{O}_{X} \ar[d]^{f_{*}\beta_{X}} \\
\tau_{Y} \ar[r]_{C(f)} & f_{*}\tau_{X}
}
\]
\end{enumerate}
\item We denote by $\cat{ISch}$
the category of idealic schemes.
\end{enumerate}
\end{Def}

\begin{Prop}
\label{prop:isch:sob:adj}
The underlying functor
$U: \cat{ISch} \to \cat{Sob}$
has a right adjoint.
\end{Prop}
\begin{proof}
We will construct a functor
$C^{+}:\cat{Sob} \to \cat{ISch}$.
Given a sober space $X$,
set $\scr{O}_{X}=\tau_{X}$,
and let $\beta_{X}$ be the identity.
Then $C^{+}(X)=(X, \tau_{X}, \Id)$ becomes
an idealic scheme.
Given a continuous map
$f:X \to Y$ of sober spaces,
set $f^{\#}=C(f):\tau_{Y} \to f_{*}\tau_{X}$.
This gives a morphism $(f,f^{\#}):C^{+}(X) \to C^{+}(Y)$
of idealic schemes.
Hence, a functor $C^{+}:\cat{Sob} \to \cat{ISch}$
is well defined.

We will show that this is the right adjoint of 
the underlying functor $U$.
The unit $\epsilon:\Id_{\cat{ISch}} \Rightarrow  C^{+} \circ U$
is defined by $\epsilon(X)=(\Id_{X},\beta_{X}\alpha)$.
The counit $U \circ C^{+} \to \Id_{\cat{Sob}}$
is given by the identity.
\end{proof}

\begin{Rmk}
Our next goal is to construct
a left adjoint of the global section functor
$\Gamma: \cat{ISch}^{\op} \to \cat{IRng$^{\dagger}$}$.
However, this seems to be impossible in general,
due to what we mentioned in the beginning
of the previous subsection.
\end{Rmk}

\subsection{$\scr{A}$-schemes}
We will introduce
a notion of $\scr{A}$-schemes,
generalizing the classical schemes
in algebraic geometry.
This is because several different
type of schemes began to appear
lately, and we thought
there should be a general theory
how to construct these schemes.

\begin{Def}
Let $\sigma$ be an algebraic type
with a multiplicative monoid structure.
\begin{enumerate}
\item A \textit{$\sigma$-algebra with a multiplicative system}
is a pair $(R,S)$, where $R$ is a $\sigma$-algebra
and $S$ is a multiplicative submonoid of $R$.
\item A \textit{homomorphism $(A,S) \to (B, T)$
of $\sigma$-algebras with multiplicative systems}
is a homomorphism $f:A \to B$ of $\sigma$-algebras 
satisfying $f(S) \subset T$.
\item We denote by $\cat{$\sigma $+Mnd}$
the category of $\sigma$-algebras with multiplicative systems.
\item For an object $(R,S) \in \cat{$\sigma $+Mnd}$,
define $L(R)=R_{S}$ as the localization of $R$ along $S$.
Note that localization exists in $\cat{$\sigma $-alg}$
for any algebraic type $\sigma$.
Given a homomorphism $f:(A,S) \to (B,T)$
of $\sigma$-algebras with multiplicative systems,
we have a natural homomorphism $L(f):A_{S} \to B_{T}$.
Hence we obtain a functor 
$L: \cat{$\sigma $+Mnd} \to \cat{$\sigma $-alg}$,
to which we refer as the \textit{localization functor}.
\end{enumerate}
\end{Def}

\begin{Def}
We will fix a quadruple $\scr{A}=(\sigma,\alpha_{1},\alpha_{2},\gamma)$
in the following. Here,
\begin{enumerate}
\item $\sigma$ is an algebraic type,
with a multiplicative monoid structure:
equivalently, there is an underlying functor
$U_{1}:\cat{$\sigma $-alg} \to \cat{Mnd}$
preserving multiplications.
\item $\alpha_{1}:\cat{$\sigma $-alg} \to \cat{PIIRng}$
is a functor, where $\cat{PIIRng}$ is the category
of pre-complete idealic semirings with idempotent multiplications.
\item $\alpha_{2}:U_{1} \Rightarrow U_{2}\alpha_{1}$
is a natural transformation,
where $U_{2}:\cat{PIIRng} \to \cat{Mnd}$
is the underlying functor, preserving multiplications.
\item The pair $\alpha=(\alpha_{1},\alpha_{2})$
gives a functor $\alpha:\cat{$\sigma $+Mnd}
\to \cat{PIIRng+Mnd}$, namely:
$\alpha(R,S)=(\alpha_{1}(R), \alpha_{2}(S))$,
and if $f:(A,S) \to (B,T)$ is a homomorphism
of $\sigma$-algebras with multiplicative systems,
then \linebreak
$(\alpha_{1}f)(\alpha_{2}(S)) \subset \alpha_{2}(T)$.

\item 
$\gamma: L \alpha_{1} \Rightarrow \alpha L$
is a natural isomorphism:
\[
\xymatrix{
\cat{$\sigma$+Mnd} \ar[r]^>(.8){\alpha} \ar[d]_{L}
\ar@{}[rd]|{\stackrel{\gamma}{\Longrightarrow}} &
\cat{PIIRng+Mnd} \ar[d]^{L} \\
\cat{$\sigma$-alg} \ar[r]_{\alpha_{1}} & \cat{PIIRng}
}
\]
\end{enumerate}
We will refer to $\scr{A}$ as a
"schematizable algebraic type".
\end{Def}

\begin{Def}
Let $\scr{A}$ be as above.
\begin{enumerate}
\item For any $\sigma$-algebra $R$
and any element $Z \in \alpha_{1}(R)$,
we denote by $R_{Z}=R_{\alpha_{2}^{-1}(Z)}$ the
\textit{localization of $R$ along $Z$}, where
\[
\alpha_{2}^{-1}(Z)=\{ x \in R \mid \alpha_{2}(x) \geq Z\}.
\]
Note that $\alpha_{2}^{-1}(Z)$ is a monoid.
\item
A schematizable algebraic type $\scr{A}$
satisfies the \textit{strong patching condition} if the following holds:

\begin{enumerate}[(i)]
\item For any $\sigma$-algebra $R$,
$\alpha_{2}(R)$ generates $\alpha_{1}(R)$ as a
pre-complete idealic semiring.
\item
Let $R$ be any $\sigma$-algebra, and
$s,s_{1},\cdots,s_{n}$ be elements
of $R$ satisfying $\alpha_{2}(s)=\sum \alpha_{2}(s_{i})$
in $\alpha_{1}(R)$.
If there are elements $a_{i} \in R_{s_{i}}$
such that $a_{i}=a_{j}$ in $R_{s_{i}s_{j}}$,
then there exists a unique
$a \in R_{s}$ such that $a=a_{i}$ in $R_{s_{i}}$.
\end{enumerate}
\item $\scr{A}$ satisfies the \textit{weak
patching condition}, if (i) holds, and (ii) holds for $s=1$.
\end{enumerate}
\end{Def}

\begin{Exam}
\label{exam:schematic:alg:type}
Here are some examples of 
schematizable algebraic types:
for any of them,
the functor $\alpha_{1}$ factors through 
$\cat{PIRng}$
(the category of pre-complete idealic semirings;
but in this subsection, we will not assume the
existence of infimum operators),
so we will just describe $\alpha_{1}^{\prime}:\cat{$\sigma $-alg}
\to \cat{PIRng}$ for (1) and (2).
\begin{enumerate}
\item The algebraic type $\sigma$ is that of rings,
and $\alpha_{1}^{\prime}:\cat{Rng} \to \cat{PIRng}$
sends a ring $R$ to the set of finitely
generated ideals on $R$.
Note that $\alpha_{1}(R)=C(\Spec R)_{\cpt}$.
A homomorphism $f:A \to B$ gives a 
homomorphism $\alpha^{\prime}_{1}(f):
\alpha^{\prime}_{1}(A) \to \alpha^{\prime}_{1}(B)$
defined by 
\[
\mathfrak{a} \mapsto f(\mathfrak{a})B.
\]
The map $\alpha_{2}:R \to \alpha_{1}(R)$
sends an element $a \in R$ to $V(a)$,
namely the support of $a$.
This gives a natural transformation
$\alpha_{2}:U_{1} \Rightarrow U_{2}\alpha_{1}$,
and a functor $\alpha=(\alpha_{1},\alpha_{2}):
\cat{Rng+Mnd} \to \cat{PIIRng+Mnd}$.

For any multiplicative system $S$ of
a ring $R$,
there is a natural isomorphism
$\gamma: \Spec R_{S} \simeq \alpha_{1}(R_{S})
\simeq \alpha_{1}(R)_{\alpha_{2}(S)}$.

In this case, the strong patching condition holds.

This schematizable algebraic type
is what we use in classical algebraic geometry.

\item The algebraic type $\sigma$ is that of monoids,
and $\alpha_{1}^{\prime}$ sends
a monoid $M$ to the set of its finitely generated ideals:
an \textit{ideal} of $M$ is a non-empty set $\mathfrak{a}$
of $M$ satisfying $ax \in \mathfrak{a}$
for any $a \in \mathfrak{a}$, $x \in M$.
The natural transformation $\alpha_{2}:M \to \alpha_{1}(M)$ 
is defined by the unit of the
adjoint $\cat{Mnd} \leftrightarrows \cat{PIRng}$, which
 sends
$a \in M$ to the saturated ideal $\overline{a}$ generated by $a$.

There is a natural isomorphism
$\gamma:\alpha_{1}(M_{S}) \simeq \alpha_{1}(M)_{\alpha_{2}(S)}$,
since $\alpha_{1}(M) =\alpha_{1}(M/M^{\times})$ for any $M$.

The weak patching condition is automatically
satisfied, since any monoid is \textit{local},
i.e. it admits a unique maximal ideal consisting
of all non-unit elements. See ~\cite{Deit}.

This schematizable algebraic type
is what we use in schemes over $\FF_{1}$.
Note that, given a ring $R$, we can construct a scheme
of this type by regarding $R$ as a multiplicative monoid.
Hence, we must fix the schematizable algebraic type to
specify how we make schemes.

\item The algebraic type $\sigma$ is 
that of pre-complete idealic semirings,
and $\alpha_{1}:\cat{PIRng} \to \cat{PIIRng}$
is the left adjoint of the underlying functor.
Note that $\alpha_{1}(R)=C(\Spec R)_{\cpt}$.
$\alpha_{2}:R \to \alpha_{1}(R)$
is the unit of the adjoint
$\cat{PIRng} \leftrightarrows \cat{PIIRng}$,
defined by $a \mapsto V(a)$.
There is a natural isomorphism
$\gamma: C(\Spec R_{S}) \simeq C(\Spec R)_{\alpha_{2}(S)}$.

The strong patching condition
is satisfied, which we already verified in
Lemma ~\ref{lem:patch:cond:idealic}.
This schematizable algebraic type is used in tropical geometry.
\end{enumerate}
\end{Exam}

\begin{Prop}
\label{prop:sheaf:spec:str}
Let $R$ be a pre-complete
idealic semiring, and $X=\Spec R$
be its spectrum.
Define a $\cat{PIRng}$-valued
presheaf $\scr{O}_{X}:C(X)_{\cpt}^{\op} \to \cat{PIRng}$
on $X^{\cpt}$ by $Z \mapsto R_{Z}$.
Then $\scr{O}_{X}$ is a sheaf.
\end{Prop}
This is a direct consequence of Lemma ~\ref{lem:patch:cond:idealic}.

\begin{Def}
Let $X$ be an algebraic sober space.
Define a $\cat{PIIRng}$-valued sheaf
$\tau_{X}^{\prime}$ on $C(X)_{\cpt}$ by
\[
Z \mapsto C(X)_{\cpt,Z}\simeq C(X\setminus Z)_{\cpt}.
\]
This is indeed a sheaf, by the above proposition.
\end{Def}

We will give a definition
of $\scr{A}$-schemes.
Note that we are defining
the structure sheaves
on $X^{\cpt}$, where $X$
is an algebraic sober space.
This is sufficient, by Proposition
~\ref{prop:equiv:sheaf:alg}.
\begin{Def}
Let $\scr{A}$ be a schematizable
algebraic type.
\begin{enumerate}
\item A triple $X=(X,\scr{O}_{X},\beta_{X})$
is a \textit{$\scr{A}$-scheme} if:
\begin{enumerate}
\item $X$ is an algebraic sober space.
\item $\scr{O}_{X}$ is a $\cat{$\sigma $-alg}$-valued
sheaf on $C(X)_{\cpt}^{\op}$.
\item $\beta_{X}:\alpha_{1}\scr{O}_{X} \to \tau^{\prime}_{X}$
is a morphism of $\cat{PIIRng}$-valued
sheaves, where $\alpha_{1}\scr{O}_{X}$
is the sheafification of $Z \mapsto \alpha_{1}(\scr{O}_{X}(Z))$.

\item The restriction map reflects localization:
Let $W \leq Z$ be any two elements of $C(X)_{\cpt}$,
and $a \in \scr{O}_{X}(Z)$ be a section
satisfying $\beta_{X} (\alpha_{2}(a)) \geq W$.
Then, the restriction map $\scr{O}_{X}(Z) \to \scr{O}_{X}(W)$
factors through $\scr{O}_{X}(Z)_{a}$:
\[
\xymatrix{
\scr{O}_{X}(Z) \ar[r]^{\rom{res}} \ar[d] &
\scr{O}_{X}(W) \\
\scr{O}_{X}(Z)_{a} \ar[ur]
}
\]
\end{enumerate}
\item A \textit{morphism 
$(X,\scr{O}_{X},\beta_{X}) \to (Y,\scr{O}_{Y},\beta_{Y})$
of $\scr{A}$-schemes}
is a pair $(f,f^{\#})$ such that:
\begin{enumerate}
\item $f:X \to Y$ is a continuous map.
\item $f^{\#}:\scr{O}_{Y} \to f_{*}\scr{O}_{X}$
is a morphism of $\cat{$\sigma $-alg}$-valued
sheaves on $Y$.
\item The following diagram of $\cat{PIIRng}$-valued
sheaves is commutative:
\[
\xymatrix{
\alpha_{1}\scr{O}_{Y} \ar[d]_{\beta_{Y}} \ar[r]^{\alpha_{1} f^{\#}} &
\alpha_{1}f_{*}\scr{O}_{X} \ar[d]^{f_{*}\beta_{X}} \\
\tau^{\prime}_{Y} \ar[r]_{C(f)_{\cpt}} & f_{*}\tau^{\prime}_{X}
}
\]
\end{enumerate}
\item We denote by $\cat{$\scr{A}$-Sch}$
the category of $\scr{A}$-schemes.
\end{enumerate}
\end{Def}

\begin{Def}
We will construct a functor
$\Spec^{\scr{A}}:\cat{$\sigma $-alg}
\to \cat{$\scr{A}$-Sch}^{\op}$ as follows:
\begin{enumerate}
\item 
Given a $\sigma$-algebra $R$,
set $X=\Spec \alpha_{1}(R)^{\dagger}$. Note that
$C(X)_{\cpt}=\alpha_{1}(R)$.
Define the structure sheaf 
$\scr{O}_{X}:\alpha_{1}(R)^{\op} \to \cat{$\sigma $-alg}$
as the sheafification of a presheaf
$\scr{O}^{\prime}_{X}:Z \mapsto R_{Z}$. 

We will define a morphism
$\beta_{X}:\alpha_{1}\scr{O}_{X} \to \tau_{X}^{\prime}$
of $\cat{PIIRng}$-valued sheaves:
For any element $Z$ of $\alpha_{1}(R)$, set 
\[
S=S(Z)=\{x \in R \mid \alpha_{2}(x) \geq Z\}.
\]
Then, we have a map
\[
\alpha_{1}\scr{O}^{\prime}_{X}(Z)
=\alpha_{1}(R_{S})
\stackrel{\gamma}{\simeq} \alpha_{1}(R)_{\alpha_{2}(S)}
\to \alpha_{1}(R)_{Z}=\tau_{X}^{\prime}(Z).
\]
Sheafifying the lefthand side, we obtain a morphism
$\beta_{X}:\alpha_{1} \scr{O}_{X} \to \tau_{X}^{\prime}$.

It is obvious that
the restriction map reflects localization.
Hence, we obtain a $\scr{A}$-scheme
$\Spec^{\scr{A}}R=(X,\scr{O}_{X},\beta_{X})$.
\item Let $\varphi:B \to A$
be a homomorphism of $\sigma$-algebras.
Set $X=\Spec^{\scr{A}} A$ and $Y=\Spec^{\scr{A}} B$.
We will construct a morphism $(f,f^{\#}):X \to Y$
as follows.
The quasi-compact continuous map $f$ is 
$\Spec(\varphi):X \to Y$.
Note that $C(f)_{\cpt}=\alpha_{1}(\varphi)$.
Let us define the morphism $f^{\#}:\scr{O}_{Y} \to f_{*}\scr{O}_{X}$
of sheaves.
For any $Z \in C(Y)$,
the homomorphism $\varphi:B \to A_{\alpha_{1}(\varphi)(Z)}$
sends any element of $S(Z)$ to an invertible element,
since $\alpha_{2}$ is a natural transformation.
Hence, it gives rise to
\[
\scr{O}^{\prime}_{Y}(Z)=B_{Z}
\to A_{\alpha_{1}(\varphi)(Z)}=\scr{O}_{X}^{\prime}(f^{-1}Z).
\]
Sheafifying the both sides,
we obtain $f^{\#}:\scr{O}_{Y} \to f_{*}\scr{O}_{X}$.
It is easy to see that the following diagram commutes:
\[
\xymatrix{
\alpha_{1}\scr{O}_{Y} \ar[r]^{\alpha_{1}f^{\#}} \ar[d]_{\beta_{Y}} &
\alpha_{1}f_{*}\scr{O}_{X} \ar[d]^{f_{*}\beta_{X}} \\
\tau_{Y}^{\prime} \ar[r]_{\alpha_{1}(\varphi)} &
f_{*}\tau_{X}^{\prime}
}
\]
Hence, we obtain a functor
$\Spec^{\scr{A}}:\cat{$\sigma $-alg} \to \cat{$\scr{A}$-Sch}$,
to which we refer to as the \textit{spectrum functor}.
\end{enumerate}
\end{Def}

\begin{Thm}
\label{thm:main:adj:spec}
Suppose $\scr{A}$ satisfies the
weak patching condition.
Then, the spectrum functor $\Spec^{\scr{A}}$ is the left adjoint of 
the global section functor
$\Gamma:\cat{$\scr{A}$-Sch}^{\op} \to \cat{$\sigma $-alg}$.
\end{Thm}
\begin{proof}
First, we will define the unit $\epsilon:\Id_{\cat{$\sigma $-alg}}
\Rightarrow  \Gamma \Spec^{\scr{A}}$.
Given a $\sigma$-algebra $R$,
set $X=\Spec^{\scr{A}} R$. Then define $\epsilon$ by
\[
R \simeq \Gamma(X, \scr{O}^{\prime}_{X}) \to \Gamma(X, \scr{O}_{X}).
\]
Next, we will define the counit
$\eta:\Spec^{\scr{A}} \Gamma \to \Id_{\cat{$\scr{A}$-Sch}^{\op}}$.
Given a $\scr{A}$-scheme $X=(X,\scr{O}_{X},\beta_{X})$,
set $Y=\Spec^{\scr{A}}\Gamma(X)$.
The continuous map $\eta:X \to \Spec^{\scr{A}} \Gamma(X)$
between the underlying spaces is induced by
$\Gamma(\beta_{X}):\alpha_{1} \Gamma(X) \to C(X)_{\cpt}$.
The morphism
$\eta^{\#}: \scr{O}_{Y} \to \eta_{*}\scr{O}_{X}$
is defined as follows:
for a given $Z \in \alpha_{1}(\Gamma(X))$,
the restriction map
$\Gamma(X) \to \scr{O}_{X}(\eta^{-1}Z)=\scr{O}_{X}(\beta_{X}(Z))$
gives rise to a homomorphism
$\scr{O}^{\prime}_{Y}(Z)=\Gamma(X)_{Z}
\to \scr{O}_{X}(\eta^{-1}Z)$,
since the restriction map reflects localization.
Sheafifying the lefthand side,
we obtain the required $\eta^{\#}$.
Thanks to the weak patching condition, all four
natural transformations
$\epsilon\Gamma$, $\Spec^{\scr{A}}\epsilon$,
$\Gamma\eta$, and $\eta\Spec^{\scr{A}}$
becomes isomorphisms.
\end{proof}

Here, we will focus on to the
case when $\scr{A}$ is induced from
idealic semirings. This case has a special feature that,
we can construct the structure
sheaf from the underlying space.
\begin{Def}
Let $\scr{A}$ be the schematizable
algebraic type induced from pre-complete
idealic semirings,
introduced in Example
~\ref{exam:schematic:alg:type}(3).
\begin{enumerate}
\item
We call $\scr{A}$-schemes 
"\textit{algebraic idealic schemes}".
We denote by $\cat{alg.ISch}$
the categories of algebraic idealic schemes.

\item We define a functor $C^{++}:\cat{alg.Sob} \to \cat{alg.ISch}$
as follows:
For an algebraic sober space $X$, set
$C^{++}(X)=(X,\tau_{X}^{\prime},\Id)$.
For a morphism $f:X \to Y$ of algebraic
sober spaces, set
$C^{++}(f)=(f,C(f)^{\cpt}):C^{++}(X) \to C^{++}(Y)$,
where $C(f)^{\cpt}:\tau^{\prime}_{Y} \to f_{*}\tau^{\prime}_{X}$
is the morphism induced by $f$.
\end{enumerate}
\end{Def}
\begin{Prop}
The functor $C^{++}$ is the right adjoint
of the underlying functor $U:\cat{alg.ISch} \to \cat{alg.Sob}$.
\end{Prop}

\begin{proof}
The unit $\epsilon:\Id_{\cat{alg.ISch}} \Rightarrow 
C^{++}U$ is given by $(\Id_{X},\beta_{X}\alpha_{1})$
for any algebraic idealic scheme $(X,\scr{O}_{X},\beta_{X})$.
The counit $\eta:UC^{++} \Rightarrow  \Id_{\cat{alg.Sob}}$
is given by the identity.
\end{proof}

\begin{Rmk}
\begin{enumerate}
\item When the schematizable algebraic type $\scr{A}$ is 
the one we introduced in Example
~\ref{exam:schematic:alg:type}(3),
then the structure sheaf $\scr{O}_{X}$
of $X=\Spec^{\scr{A}} R$ is the
functor defined by $Z \mapsto R_{[Z]}$ for $Z \in C(X)$.
This follows from Lemma ~\ref{lem:patch:cond:idealic}.
\item
Note that we don't have a natural
underlying functor $\cat{alg.ISch} \to \cat{ISch}$:
The sheaves $\tau$ and $\tau^{\prime}$ which
represent the topological structure, are different.
\end{enumerate}
\end{Rmk}

Let us summarize all the categories and functors
we have obtained so far:
Let $\scr{A}=(\sigma,\alpha_{1},\alpha_{2},\gamma)$
be a schematizable algebraic type.
Then the functors are illustrated below:
The pairs of arrows are adjoints.
The equal signs are equivalences.
\[
\xymatrix{
&&& \cat{alg.IRng$^{\dagger}$}
\ar[r]^{\Spec} \ar[d] &
\cat{alg.ISch}^{\op} \ar@<.5ex>[l]^{\Gamma}
\ar[d]^{U} \\
\cat{alg.$\sigma^{\dagger}$-alg} 
\ar@{=}[r]_{\comp^{\prime}}^{U_{\cpt}}
\ar[d]^{U} &
\cat{$\sigma$-alg} \ar@<.5ex>[dl]^{\comp}  
\ar[r]^{\Spec^{\scr{A}}} &
\cat{$\scr{A}$-Sch}^{\op} \ar@<.5ex>[l]^{\Gamma}
\ar[r]^{U} \ar@{.>}[ur]^{U^{\prime}}&
\cat{alg.IIRng$^{\dagger}$} \ar@<.5ex>[u]
\ar@{=}[r]^{\Spec}_{C} 
\ar[d]^{U} &
\cat{alg.Sob}^{\op} \ar@<.5ex>[u]^{C^{++}}
\ar[d]^{U} \\
\cat{$\sigma^{\dagger}$-alg} \ar@<.5ex>[u]^{\alg}
\ar[ur]^>(.45){U}
&&&
\cat{IIRng$^{\dagger}$} \ar@<.5ex>[u]^{\alg} \ar[r]^{\Spec} &
\cat{Sob}^{\op} \ar@<.5ex>[u]^{\alg} \ar@<.5ex>[dl]^{C^{+}}
\ar@<.5ex>[l]^{C} \ar[d]^{U} \\
&&& \cat{ISch}^{\op} \ar[ur]^{U}
&\cat{Top}^{\op} \ar@<.5ex>[u]^{\sob}
}
\]
As we see, the right half can be
regarded as the categories of geometric objects,
while the left half are those of algebraic objects.
The most mysterious part is the functor
$U:\cat{$\scr{A}$-Sch} \to \cat{alg.IIRng$^{\dagger}$}$
in the middle, and it has been (and will continue to be)
the central subject of algebraic geometry,
arithmetic geometry, and tropical geometry.
If $\scr{A}$ is one of the example in 
\ref{exam:schematic:alg:type},
then $U$ factor through $\cat{alg.IRng$^{\dagger}$}$.

\subsection{Comparison with classical schemes}
Next, we will see the relation between the
category of classical schemes and
the category of new schemes which we have introduced.
In the proceedings,
$\scr{A}$ is the schematizable
algebraic type induced from rings,
introduced in Example
~\ref{exam:schematic:alg:type}(1).

\begin{Prop}
Let $\cat{QSch}$ be the category of 
quasi-compact, quasi-seperated schemes
and quasi-compact morphisms, in the
classical sense.
Then, there is a fully faithful functor
$I:\cat{QSch} \to \cat{$\scr{A}$-Sch}$
which satisfies $I \Spec^{\op} \simeq \Spec^{\scr{A},\op}$.
\[
\xymatrix{
\cat{Rng}^{\op} \ar[r]^{\Spec^{\op}} \ar[rd]_{\Spec^{\scr{A},\op}} &
\cat{QSch} \ar[d]^{I} \\
& \cat{$\scr{A}$-Sch}
}
\]
\end{Prop}
\begin{proof}
Let $X=(X,\scr{O}_{X})$ be a quasi-compact, quasi-seperated
scheme in the classical sense.
We must construct a morphism 
$\beta_{X}:\alpha_{1}\scr{O}_{X} \to \tau_{X}^{\prime}$
of sheaves.
Since $X$ is locally isomorphic to an affine scheme,
we already have a local isomorphism
$\alpha_{1}\scr{O}_{X} \simeq \tau_{X}^{\prime}$,
which patch up to give a global morphism $\beta_{X}$.
Hence we obtain a $\scr{A}$-scheme $(X, \scr{O}_{X},\beta_{X})$.

For any quasi-compact morphism $f:X \to Y$ of quasi-compact,
quasi-seperated
schemes, $f$ commutes with $\beta$,
since $f$ is locally induced by a homomorphism
of rings. This means that
$f$ naturally becomes a morphism
of $\scr{A}$-schemes.
Hence, we have a functor $I:\cat{QSch} \to \cat{$\scr{A}$-Sch}$.

We claim that $I$ is fully faithful.
Let $X, Y$ be quasi-compact, quasi-seperated schemes,
and $f:I(X) \to I(Y)$ be a morphism
of $\scr{A}$-schemes.
We already know that
$I(X)$ and $I(Y)$ are locally ringed spaces,
and it suffices to show that
$f$ is a morphism of locally ringed spaces,
i.e. $\scr{O}_{Y,f(x)} \to \scr{O}_{X,x}$
is a local homomorphism for any $x \in X$.
This is a local argument, hence
we may assume $X$ and $Y$ are both affine,
say $X=\Spec A$ and $Y=\Spec B$.
Using Theorem \ref {thm:main:adj:spec}, we have
\[
\Hom_{\cat{$\scr{A}$-Sch}}(I(X),I(Y))
\simeq \Hom_{\cat{Rng}}(B,A)
\simeq \Hom_{\cat{QSch}}(X,Y),
\]
which shows that $f$ is induced locally
(and hence globally)
from a morphism of quasi-compact, quasi-seperated schemes.

It is clear from the construction of $\Spec^{\scr{A}}$
that $I\Spec^{\op}\simeq \Spec^{\scr{A},\op}$.
\end{proof}

\begin{Rmk}
Let $\cat{$\scr{M}$-Sch}$ be the
category of quasi-compact, quasi-seperated
schemes over $\FF_{1}$,
in the sense of Deitmar \cite{Deit}.
and $\scr{A}$ be the schematizable algebraic type
of Example ~\ref{exam:schematic:alg:type}(2).
Then,
by the same arguments as above,
we obtain a fully faithful functor
$\cat{$\scr{M}$-Sch} \to \cat{$\scr{A}$-Sch}$.
\end{Rmk}

\textsc{S. Takagi: Department of Mathematics, Faculty of Science,
Kyoto University, Kyoto, 606-8502, Japan}

\textit{E-mail address}: \texttt{takagi@math.kyoto-u.ac.jp}
\end{document}